\documentclass{article}
\usepackage{amsfonts}                     
\usepackage{times}                        
\usepackage{fancyheadings}                
\usepackage{amsmath, amsthm, amssymb} 
\usepackage[sc]{mathpazo} 
\usepackage[T1]{fontenc} 
\usepackage{microtype} 

\usepackage[hmarginratio=1:1,top=32mm,columnsep=20pt]{geometry} 
\usepackage[hang, small,labelfont=bf,up,textfont=rm,up]{caption} 
\usepackage{booktabs} 
\usepackage{float} 
\usepackage{hyperref} 
\usepackage{epsfig} 
\usepackage{epstopdf} 
\usepackage{multirow}
\usepackage{tikz}	
\usetikzlibrary{arrows}
\usepackage{verbatim}
\usepackage[parfill]{parskip}
\usepackage[symbol]{footmisc}	
\usepackage{enumitem}
\numberwithin{equation}{section}

\usepackage{tikz}
\usetikzlibrary{arrows}

\newcommand{\beq}{\begin{equation*}}
\newcommand{\eeq}{\end{equation*}}
\newcommand{\beqa}{\begin{align*}}
\newcommand{\eeqa}{\end{align*}}
\newcommand{\ben}{\begin{enumerate}}
\newcommand{\een}{\end{enumerate}}
\newcommand{\bed}{\begin{definition}}
\newcommand{\eed}{\end{definition}}
\newcommand{\bet}{\begin{theorem}}
\newcommand{\eet}{\end{theorem}}
\newcommand{\bel}{\begin{lemma}}
\newcommand{\eel}{\end{lemma}}
\newcommand{\bec}{\begin{corollary}}
\newcommand{\eec}{\end{corollary}}
\newcommand{\bep}{\begin{proof}}
\newcommand{\eep}{\end{proof}}

\newcommand{\tab}{\hspace*{2em}}
\newcommand{\tbl}{\textquotedblleft}
\newcommand{\tbr}{\textquotedblright}
\newcommand{\modu}[1]{\ensuremath{(\!\!\!\!\mod #1)}}
\newcommand{\floor}[1]{\ensuremath{\left\lfloor#1\right\rfloor}}
\newcommand{\ceil}[1]{\ensuremath{\left\lceil#1\right\rceil}}

\newcommand{\seqnum}[1]{\href{http://oeis.org/#1}{\underline{#1}}}

\def\Z{{\mathbb Z}}

\def\N{{\mathbb{N}}}
\newtheorem{theorem}{Theorem}

\newtheorem{conjecture}[theorem]{Conjecture}
\newtheorem{corollary}[theorem]{Corollary}

\newtheorem{lemma}[theorem]{Lemma}

\newtheorem{proposition}[theorem]{Proposition}

\theoremstyle{definition}
\newtheorem{definition}[theorem]{Definition}
\newtheorem{example}[theorem]{Example}

\begin{document}                          

\title{$q$-Binomials and related symmetric unimodal polynomials}

\author{Bryan Ek\footnote{Department of Mathematics, The School of Arts and Sciences, Rutgers, The State University of New Jersey, Piscataway, NJ 08854}}
\renewcommand{\thefootnote}{\arabic{footnote}}
\setcounter{footnote}{0}

\maketitle

\begin{abstract} \noindent
	The $q$-binomial coefficients were conjectured to be unimodal as early as the 1850's, but it remained unproven until Sylvester's 1878 proof using invariant theory. In 1982, Proctor gave an \tbl elementary\tbr\ proof using linear algebra. Finally, in 1989, Kathy O'Hara provided a {\it combinatorial} proof of the unimodality of the $q$-binomial coefficients. Very soon thereafter, Doron Zeilberger translated the argument into an elegant recurrence. We introduce several perturbations to the recurrence to create a larger family of unimodal polynomials. We analyze how these perturbations affect the final polynomial and analyze some specific cases.
\end{abstract}


\begin{section}{Introduction}

	\begin{subsection}{Motivation}
		\begin{quote}
			\tbl The study of unimodality and log-concavity arise often in combinatorics, economics of uncertainty and information, and algebra, and have been the subject of considerable research.\tbr \cite{Alvarez2000}
		\end{quote}
		Intuitively, many sequences seem to be unimodal, but how does one prove that fact? We will review some methods of building unimodal sequences as well as a few lemmas that imply unimodality.
		
		Knowing a sequence is unimodal allows for guaranteed discovery of the global extremum using an easy search algorithm. Unimodality is also useful for probability applications. Identifying a probability distribution as unimodal allows certain approximations for how far a value will be from its mode (Gauss' inequality \cite{GaussInequality}) or mean (Vysochanskij-Petunin inequality \cite{VPinequality}).\\
		
		Suppose we start with a set of nice combinatorial objects that satisfy property $\mathcal{P}$. How can we make the set larger and still satisfy property $\mathcal{P}$? Or a slightly different property $\mathcal{P}'$? The goal of this project is to use reverse engineering to obtain highly non-trivial, surprising theorems about unimodality.\\

		\begin{example}
			Consider the following functions
			\begin{align*}
				P_1(n)	&=	{\frac {1-{q}^{n+1}}{1-q}},\\
				P_2(n)	&=	{\frac {1- q^{2n+1}}{1-q}},\\
				P_3(n)	&=	{\frac {2q^{3n+2}-2q^{3n+1}+q^{3n}- q^{2n+1}-q^{n+1}+{q}^{2}-2\,q+2}{ \left( 1-q \right) ^{2}}},\\
				P_4(n)	&=	\bigg(5q^{4n+2}-5q^{4n+1}+3{q}^{4n-2}-3q^{4n-3}+4q^{4n-4}-4{q}^{2n+3}-4{q}^{2n-1}\\
						&\tab	+4{q}^{6}-3{q}^{5}+3{q}^{4}-5q+5\bigg)\bigg/\left(1-q\right) ^{2},\\
				P_5(n)	&=	\bigg(2q^{5n+3}-4q^{5n+2}+7{q}^{5n+1}-5{q}^{5n}+5{q}^{5n-3}-5{q}^{5n-4}+4{q}^{5n-5}-4{q}^{5n-6}+6{q}^{5n-9}\\
						&\tab	-6q^{5n-10}+3q^{5n-11}	-5{q}^{4n+2}+5{q}^{4n+1}-4{q}^{4n-2}+4{q}^{4n-3}-3{q}^{4n-4}-5{q}^{3n+4}\\
						&\tab	+5{q}^{3n+3}-3{q}^{3n+2}-6{q}^{3n-4}+6{q}^{3n-5}-3{q}^{3n-6}+3{q}^{2n+9}-6{q}^{2n+8}+6{q}^{2n+7}\\
						&\tab+3{q}^{2n+1}-5{q}^{2n}+5{q}^{2n-1}+3{q}^{n+7}-4{q}^{n+6}+4{q}^{n+5}-5{q}^{n+2}+5{q}^{n+1}\\
						&\tab	-3{q}^{14}+6{q}^{13}-6{q}^{12}+4{q}^{9}-4{q}^{8}+5{q}^{7}-5{q}^{6}+5{q}^{3}-7q^{2}+4q-2\bigg)\bigg/\left( -1+q \right) ^{3}.
			\end{align*}
			All of the preceding functions are not only polynomial for $n\ge0,0,0,4,4$, respectively, but they are also unimodal.
			\label{exa:ImpressivePolynomials1}
		\end{example}
		
		\begin{example}
			The functions
			\begin{align*}
				Q_1(n)	&=	{\frac {4\left(1-{q}^{n+1}\right)}{1-q}},\\
				Q_2(n)	&=	5{\frac {8{q}^{2n+4}-2{q}^{n+3}-2{q}^{n+2}+8q-3[(q^5+1)(q^n-(-q)^n)+(q^4+q)(q^n+(-q)^n)]}{2q\left(1-q\right)\left(1-q^2\right)}},\\
				Q_3(n)	&=	\bigg(16\big[3-3q+16{q}^{2}-16{q}^{n+1}-16{q}^{n+3}-16{q}^{n+5}+16{q}^{2n+3}\\
						&\tab\tab	+16{q}^{2n+5}+16{q}^{2n+7}-16{q}^{3n+6}+3{q}^{3n+7}-3{q}^{3n+8}\big]\\
						&\tab+(1-(-1)^n)(-q)^{(3n-9)/2}\big[8{q}^{17}-8{q}^{16}+64{q}^{15}-9{q}^{14}-55{q}^{13}-60{q}^{11}+57{q}^{10}\\
						&\tab\hspace{3cm}	-125{q}^{9}+125{q}^{8}-57{q}^{7}+60{q}^{6}+55{q}^{4}+9{q}^{3}-64{q}^{2}+8{q}-8\big]\\
						&\tab+(1-(-1)^n)q^{(3n-9)/2}\big[8{q}^{17}-8{q}^{16}+64{q}^{15}+9{q}^{14}-73{q}^{13}-60{q}^{11}-65{q}^{10}\\
						&\tab\hspace{3.7cm}	-3{q}^{9}+3{q}^{8}+65{q}^{7}+60{q}^{6}+73{q}^{4}-9{q}^{3}-64{q}^{2}+8q-8\big]\\
						&\tab+(1+(-1)^n)(-q)^{(3n-6)/2}\big[12{q}^{14}-12{q}^{13}+64{q}^{12}-75{q}^{11}+74{q}^{10}-127{q}^{9}\\
						&\tab\hspace{4.4cm}	+127{q}^{5}-74{q}^{4}+75{q}^{3}-64{q}^{2}+12q-12\big]\\
						&\tab+(1+(-1)^n)q^{(3n-6)/2}\big[12{q}^{14}-12{q}^{13}+64{q}^{12}-53{q}^{11}-74{q}^{10}-{q}^{9}\\
						&\tab\hspace{3.2cm}+{q}^{5}+74{q}^{4}+53{q}^{3}-64{q}^{2}+12q-12\big]\bigg)
						\bigg/4(1-q)^2(1-q^6),
			\end{align*}
			are all unimodal for $n\ge0$. 
			\label{exa:ImpressivePolynomials2}
			The function $Q_3$ is truly amazing. It appears quite unwieldy and at first glance one might doubt that it even has real coefficients let alone integer coefficients. But for any nonnegative integer $n$, $Q_3(n)$ is guaranteed to be a unimodal polynomial in $\Z[q]$. Try simplifying $Q_3$ assuming $n$ is even or odd. One can obtain further simplifications assuming $n=0,1,2,3\modu{4}$.
		\end{example}

	\end{subsection}

	\begin{subsection}{Symmetric and Unimodal}
		We recall several definitions and propositions from Zeilberger \cite{Zeilberger1989} for the sake of completeness.
		\bed[Unimodal]
			A sequence $\mathcal{A}=\{a_0,\ldots,a_n\}$ is {\it unimodal} if it is weakly increasing up to a point and then weakly decreasing, i.e., there exists an index $i$ such that $a_0\le a_1\le\cdots\le a_i\ge\cdots\ge a_n$.
		\eed
		\bed[Symmetric]
			A sequence $\mathcal{A}=\{a_0,\ldots,a_n\}$ is {\it symmetric} if $a_i=a_{n-i}$ for every $0\le i\le n$.
		\eed
		A polynomial is said to have either of the above properties if its sequence of coefficients has the respective property.
		
		\bed[Darga]
			The {\it darga} of a polynomial $p(q)=a_iq^i+\cdots+a_jq^j$, with $a_i\ne0\ne a_j$, is defined to be $i+j$, i.e., the sum of its lowest and highest powers.
		\eed
		Brent \cite{BRENTI1990} uses $C(p)$, that is the average of lowest and highest powers. I.e., $C(p)=\frac{darga(p)}{2}$.
		\begin{example}
			$darga(q^2+3q^3)=5$ and $darga(q^2)=4$.
		\end{example}
		
		
		\begin{proposition}
			The sum of two symmetric and unimodal polynomials of darga $m$ is also symmetric and unimodal of darga $m$.
			\label{pro:SumSymUni}
		\end{proposition}
		\begin{proposition}
			The product of two symmetric and unimodal nonnegative\footnote{Nonnegative is necessary: $(-1+x-x^2)^2=1-2x+3x^2-2x^3+x^4$.} polynomials of darga $m$ and $m'$ is a symmetric and unimodal polynomial of darga $m+m'$.
			\bep
				A polynomial is symmetric and unimodal of darga $m$ if and only if it can be expressed as a sum of \tbl atomic\tbr entities of the form $c(q^{m-r}+q^{m-r-1}+\cdots+q^{r})$, for some {\bf positive} constant $c$ and integer $0\le r\le\frac{m}{2}$. By Proposition \ref{pro:SumSymUni}, it is enough to prove that the product of two such atoms of dargas $m$ and $m'$ is symmetric and unimodal of darga $m+m'$;
				\beq
					(q^{m-r}+\cdots+q^{r})(q^{m'-r'}+\cdots+q^{r'})=q^{m+m'-r-r'}+2q^{m+m'-r-r'-1}+\cdots+2q^{r+r'+1}+q^{r+r'}.
				\eeq
			\eep
			\label{pro:ProdSymUni}
		\end{proposition}
		\begin{proposition}
			If $p$ is symmetric and unimodal of darga $m$, then $q^\alpha p$ is symmetric and unimodal of darga $m+2\alpha$.
			\label{pro:PowerSymUni}
		\end{proposition}
		One example is the binomial polynomial: $(1+x)^m$. It is symmetric and unimodal of darga $m$.
		
		\bed[$\gamma$-nonnegative \cite{Branden2015}]
			If $h(x)$ is a symmetric function (in its coefficients), then we can write $h(x)=\sum_{k=0}^{\floor{m/2}}\gamma_kx^k(1+x)^{m-2k}$. We call $\{\gamma_k\}_{k=0}^{\floor{m/2}}$ the {$\gamma$-vector} of $h$. If the $\gamma$-vector is nonnegative, $h$ is said to be {\it $\gamma$-nonnegative}.
		\eed
		One can use Propositions \ref{pro:SumSymUni}, \ref{pro:ProdSymUni}, and \ref{pro:PowerSymUni} to prove that $\gamma$-nonnegative implies symmetric and unimodal of darga $m$.
		
		\begin{example}
			To prove the polynomials in Example \ref{exa:ImpressivePolynomials1} are actually unimodal, one can show
			\begin{align*}
				P_1(n)	&=	\sum_{i=0}^{n}{q}^{i},\\
				P_2(n)	&=	P_1(2n),\\
				P_3(n)	&=	2P_1\left(3n\right) +{q}^{2}P_1\left(2n-2\right)P_1\left(n-2\right),\\
				P_4(n)	&=	5P_1\left(4n\right)+3{q}^{4}P_2\left(2n-4\right)+4{q}^{6}P_1\left(2n-4\right)P_2\left(n-4\right),\\
				P_5(n)	&=	2P_1\left(5n\right)+5{q}^{2}P_1\left(4n-2\right)P_1\left(n-2\right)+4{q}^{8}P_2\left(2n-6\right)P_1\left(n-4\right)\\
						&\tab	+5{q}^{6}P_1\left(3n-4\right)P_2\left(n-4\right)+3{q}^{12}P_1\left(2n-6\right)P_3\left(n-6\right),
			\end{align*}
			and then use Propositions \ref{pro:SumSymUni}, \ref{pro:ProdSymUni}, and \ref{pro:PowerSymUni} along with the property that $darga(P_k(n))=nk$. One must initially assume $n\ge 0,0,2,4,6$, respectively, but the $k=3,5$ bounds can be lowered by checking smaller values of $n$ manually.
		\end{example}
		
		\begin{example}
			To prove the polynomials in Example \ref{exa:ImpressivePolynomials2} are in fact unimodal, one can show
			\begin{align*}
				Q_1(n)	&=	4\sum_{i=0}^{n}{q}^{i},\\
				Q_2(n)	&=	{q}^{2}Q_2\left(n-2\right)+5Q_1\left(2n\right),\\
				Q_3(n)	&=	{q}^{6}Q_3\left(n-4\right)+3Q_1\left(3n\right)+4{q}^{2}Q_1\left(2n-2\right)Q_1\left(n-2\right),
			\end{align*}
			and then use induction on $n$ for $Q_2,Q_3$, with base cases $n=\{0,1\},\{0,1,2,3\}$, respectively, and Propositions \ref{pro:SumSymUni}, \ref{pro:ProdSymUni}, and \ref{pro:PowerSymUni} along with the property that $darga(Q_k(n))=nk$.
		\end{example}
		The wonderful aspect about these proofs is that they can be easily verified by computer!
		
		\smallbreak
		
		Being able to write a function in \tbl decomposed\tbr\ form allows us to quickly verify the unimodal nature of each part and hence the sum. However, the combined form generally does not lead to an obvious decomposition. This is the motivation behind using reverse engineering: to guarantee that the decomposed form exists.
		
		\begin{subsubsection}{Real-rootedness and Log-concavity}
			The following properties are greatly related to unimodality but turn out to not be applicable in our situation. They are included for completeness of discussion.
			\bed[Real-rootedness]
				The generating polynomial, $p_\mathcal{A}(x) := a_0 + a_1q +\cdots+ a_nq^n$, is called {\it real-rooted} if all its zeros are real. By convention, constant polynomials are considered to be real-rooted.
			\eed
			\bed[Log-concavity]
				A sequence $\mathcal{A}=\{a_0,\ldots,a_n\}$ is {\it log-concave (convex)} if $a_i^2\ge a_{i-1}a_{i+1}$ ($a_i^2\le a_{i-1}a_{i+1}$) for all $1\le i<n$.
			\eed
			There is also a notion of {\it $k$-fold log-concave} and {\it infinitely log-concave} with open problems that may be of interest to the reader \cite{Branden2015}.
			\begin{proposition}
				The Hadamard (term-wise) product of log-concave (convex) sequences is also log-concave (convex).
			\end{proposition}
			\bel[Br\"{a}nd\'en \cite{Branden2015}]
				Let $\mathcal{A}=\{a_k\}_{k=0}^n$ be a finite sequence of nonnegative numbers.
				\begin{itemize}
					\item	If $p_\mathcal{A}(x)$ is real-rooted, then the sequence $\mathcal{A}':=\{a_k/\binom{n}{k}\}_{k=0}^n$ is log-concave.
					\item	If $\mathcal{A}'$ is log-concave, then so is $\mathcal{A}$.
					\item	If $\mathcal{A}$ is log-concave and positive,\footnote{Nonnegative is not sufficient: $\{1,0,0,1\}$ is log-concave (and log-convex) but not unimodal.} then $\mathcal{A}$ is unimodal.
				\end{itemize}
				\label{lem:Branden}
			\eel
			For a self-contained proof using less general (and possibly easier to understand) results, see Lecture 1 from Vatter's Algebraic Combinatorics class \cite{Vatter2009}, which uses the book {\it A Walk Through Combinatorics} \cite{Bona2006}. The converse statements of Lemma \ref{lem:Branden} are false. There are log-concave polynomials that are not real-rooted and there are unimodal polynomials which are not log-concave.
			
			Br\"and\'en also references a proof by Stanley \cite{Stanley1989} that:
			\begin{lemma}
				If $A(x),B(x)$ are log-concave, then $A(x)B(x)$ is log-concave. And if $A(x)$ is log-concave and $B(x)$ is unimodal, then $A(x)B(x)$ is unimodal.
				\label{lem:Stanley}
			\end{lemma}
			It is not sufficient for $A,B$ to be unimodal: $(3+q+q^2)^2=9+6q+7q^2+2q^3+q^4$.\footnote{Stanley \cite{Stanley1989} incorrectly writes the coefficient of $q^3$ as $1$.}
		\end{subsubsection}
		
	\end{subsection}

	\begin{subsection}{Partitions}
		\bed[Partition]
			A {\it partition} of $k$ is a non-increasing sequence of positive integers $\lambda=[a_1,a_2,\ldots,a_s]$ s.t. $\sum_{i=1}^s a_i=k$.\\
			We will use {\it frequency representation} $\lambda=[a_1^{b_1},a_2^{b_2},\ldots,a_r^{b_r}]$ s.t. $a_i>a_{i+1}$ and $b_i>0$ to abbreviate repeated terms \cite{Andrews1998}.\footnote{If $b_i=1$, it is omitted.} Note now that $\sum_{i=1}^r a_i b_i=k$. There is possible ambiguity as to whether $x^y$ is $x$ repeated $y$ times, or $x^y$ counted once. We always reserve exponentials in partitions to indicate repetition.\\
			The {\it size} of a partition, denoted $|\lambda|$, will indicate the {\bf number} of parts.\footnote{Many papers use $|\lambda|=k$ to denote what $\lambda$ partitions.} In standard notation $|\lambda|=s$; in frequency notation $|\lambda|=\sum_{i=1}^r b_i$.\\
			The number of parts of size $i$ in the partition is denoted $d_i$. In frequency notation, $d_{a_i}=b_i$.\\
			To indicate that $\lambda$ is a partition of $k$, we use $\lambda\vdash k$.
		\eed
		We will use $\lambda$ to denote a partition. Unless otherwise stated, $\lambda$ is a partition of $k$.
	\end{subsection}

	\begin{subsection}{Paper Organization}
		This paper is organized in the following sections:
		\ben[label=\bfseries Section \arabic*:]
			\setcounter{enumi}{1}
			\item Maple Program: Briefly describes the accompanying Maple package.
			\item	$q$-binomial Polynomials: Describes the original interesting polynomials.
			\item	Original Recurrence: Introduces the recurrence that generates the $q$-binomial polynomials and simultaneously proves their unimodality. We briefly discuss properties of the recurrence itself.
			\item	Altered Recurrence: Modifications are injected into the recurrence that still maintain unimodality for the resulting polynomials. We examine the effects various changes have. This section contains most of the reference to new unimodal polynomials.
			\item	OEIS: Uses some of the created unimodal polynomials to enumerate sequences for the OEIS.\footnote{Online Encyclopedia of Integer Sequences \cite{OEIS}.}
			\item	Conclusion and Future Work: Provides avenues for future research.
		\een
	\end{subsection}
	
\end{section}

\begin{section}{Maple Program}
	The backbone of this paper is based on experimental work with the Maple package {\bf Gnk} available at
	\begin{center}{\bf math.rutgers.edu/$\sim$bte14/Code/Gnk/Gnk.txt}.\end{center}
	I will mention how functions are formatted throughout this paper using {\sc implemented as {\bf function}}. For general package help and a list of available functions, type {\bf Help()}. For help with a specific function, type {\bf Help(function)}.
	
	The most important function is {\bf KOHgeneral}. It modifies the original recurrence in Eqn.\ \eqref{eqn:KOH} in several different ways. Using this function, one can change the recurrence call, multiply the summand by any manner of constant, restrict the types of partitions, or even hand-pick which partitions should be weighted most highly. The key part of these changes is that the solution to the recurrence remains symmetric and unimodal of darga $nk$.
	
	{\bf math.rutgers.edu/$\sim$bte14/Code/Gnk/Polynomials/} contains examples of many non-trivial one-variable rational functions that are guaranteed to be unimodal polynomials (for $n\ge n_0$). They were created using the {\bf KOHrecurse} function and restricting the partitions summed over in \eqref{eqn:KOH} to have smallest part $\ge1,2,3,4$ and distance $\ge1,2,3,4$. In decomposed form, they are clearly unimodal; but in combined form one is hard-pressed to state unimodality with certainty.
	
	The \tbl impressive\tbr\ polynomials in Examples \ref{exa:ImpressivePolynomials1} and \ref{exa:ImpressivePolynomials2} were generated by {\bf RandomTheoremAndProof}. It uses the recurrence in Eqn.\ \eqref{eqn:KOH} with random multiplicative constants in each term. The parameters were $(n,q,5,5,0,false)$ and $(n,q,3,5,0,true)$, respectively. Initially, $Q_3(n)$ was an even larger behemoth. This will typically happen when Maple solves a recurrence equation so going beyond $K=3$ and setting the argument {\em complicated} to {\em true} is not recommended unless one wants to spend substantial time parsing the polynomial into something much more readable. Or better yet, one could come up with a way for Maple to do this parsing automatically.
	
	There are functions to test whether a polynomial is symmetric, unimodal, or both ({\bf isSymmetric(P,q)}, {\bf isUnimodal(P,q)}, {\bf isSymUni(P,q)}, respectively). As an extra take-away, {\bf generalPartitions} outputs all partitions restricted to a minimum/maximum integer/size as well as difference\footnote{Consecutive parts in the partition differ by $\ge d$.} and congruence requirements. One can also specify a finite set of integers which are allowed in {\bf genPartitions}.
\end{section}

\begin{section}{$q$-binomial Polynomials}
	Of particular interest among unimodal polynomials are the $q$-binomial polynomials ${{n}\brack{k}}_q$ {\sc implemented as {\bf qbin(n,k,q)}}. The $q$-binomials are analogs of the binomial coefficients. In the limit $q\to1^-$, we obtain the usual binomial coefficients. We parametrize $q$-binomial coefficients in order to eliminate the redundancy: ${{n}\brack{k}}_q={{n}\brack{n-k}}_q$. The function of interest is
	\begin{equation}
		G(n,k)	=	G(k,n)	=	{{n+k}\brack{k}}_q	=	\frac{[n+k]_q!}{[n]_q![k]_q!}	=	\frac{(1-q^{n+1})\cdots(1-q^{n+k})}{(1-q)\cdots(1-q^k)},
		\label{eqn:Qbinomial}
	\end{equation}
	and is {\sc implemented as {\bf Gnk(n,k,q)}}. $[n]_q=\frac{1-q^n}{1-q}$ is the q-bracket. $[n]_q!=\prod_{i=1}^n[i]_q$ denotes the q-factorial\footnote{The q-factorial is the generating polynomial for the number of inversions over the symmetric group $S_n$ \cite{Branden2015}.} {\sc implemented as {\bf qfac(n,q)}}. It is non-trivial that $G(n,k)$ is even a polynomial for integer $n,k$. However, this follows since $G(n,k)$ satisfies the simple recurrence relation
	\begin{equation}
		G(n+1,k+1)	=	q^{k+1}G(n,k+1)+G(n+1,k).
		\label{eqn:Recurrence}
	\end{equation}
	This recurrence also shows that $G(n,k)$ is the generating function for the number of partitions in an $n\times k$ box. A partition either has a part of largest possible size ($q^{k+1}G(n,k+1)$) or it does not ($G(n+1,k)$). To prove this, one can simply check that $G(n,k)$ does satisfy this recurrence as well as the initial conditions $G(n,0)=1=G(0,k)$.\footnote{The coefficients are symmetric since each partition of $\ell$ inside an $n\times k$ box corresponds to a partition of $nk-\ell$.} This argument also confirms $G(n,k)=G(k,n)$.\footnote{$G(n-k,k)$ also counts the number of subspaces of dimension $k$ in a vector space of dimension $n$ over a finite field with $q$ elements.} 

	The question of unimodality for the $q$-binomial coefficients was first stated in the 1850s by Cayley and then proven by Sylvester in 1878 \cite{Branden2015}. The first \tbl elementary\tbr\ proof was given by Proctor in 1982. He essentially described the coefficients as partitions in a box (grid-shading problem) and then used linear algebra to finish the proof \cite{Proctor1982}.
	
	Lemma \ref{lem:Branden} does not apply here; q-brackets and q-factorials are in the class of log-concave (use Lemma \ref{lem:Stanley}) but not real-rooted (for $n>1$) polynomials while these $G(n,k)$ polynomials are in the class of unimodal, but not necessarily log-concave polynomials. For example, $G(2,2)=1+q+2q^2+q^3+q^4$ is unimodal but not log-concave. In fact, most $q$-binomial polynomials are not log-concave since the first (and last) 2 coefficients are 1. Explicitly: $G(n,k)$ is log-concave if and only if $k=1$.
	
\end{section}

\begin{section}{Original Recurrence}
	The symmetric and structured nature of $G(n,k)$ and its coefficients is visually muddled by the following recurrence that Zeilberger \cite{Zeilberger1989} created to translate O'Hara's \cite{Ohara} combinatorial argument\footnote{Her argument is rewritten and explained by Zeilberger \cite{zeilberger1989kathy}.} into mostly algebra as
	\begin{equation}
		G(n,k)	=	\sum_{(d_1,\ldots,d_k);\sum_{i=1}^k id_i=k} q^{k\left(\sum_{i=1}^kd_i\right)-k-\sum_{1\le j<i\le k}(i-j)d_id_j}	\prod_{i=0}^{k-1}G\left((k-i)n-2i+2\sum_{j=0}^{i-1}(i-j)d_{k-j},d_{k-i}\right)
		\label{eqn:KOH}
	\end{equation}
	is {\sc implemented as {\bf KOH(n,k,q)}}. The product is part of the summand and the outside sum is over all partitions of $k$. $d_i$ is the number of parts of size $i$ in the partition. The initial conditions\footnote{$G(n,1)$ is explicitly given because Eqn.\ \eqref{eqn:KOH} only yields the tautology $G(n,1)=G(n,1)$.} are
	\begin{equation}
		G(n<0,k)	=	0,	\hspace{1.35em}	G(n,k<0)	=	0,	\hspace{1.35em}	G(0,k)	=	1,	\hspace{1.35em}	G(n,0)	=	1,	\hspace{1.35em}	G(n,1)	=	\frac{1-q^{n+1}}{1-q}.
		\label{eqn:InitialConditions}
	\end{equation}
	By using Propositions \ref{pro:SumSymUni}, \ref{pro:ProdSymUni}, and \ref{pro:PowerSymUni}, and induction on the symmetric unimodality of $G(a,b)$ for $a<n$ or $b<k$, we can see (after a straightforward calculation) that the right hand side is symmetric and unimodal of darga $nk$. For each partition, the darga will be
	\beq
		2\left[k\left(\sum_{i=1}^kd_i\right)-k-\sum_{1\le j<i\le k}\hspace{-0.8em}(i-j)d_id_j\right]+\sum_{i=0}^{k-1}\left((k-i)n-2i+2\sum_{j=0}^{i-1}(i-j)d_{k-j}\right)d_{k-i}	=	nk.
	\eeq
	Recall that $\sum_{i=1}^k id_i=k$.
	
	\bigskip
	
	One is led to ask how long this more complicated recurrence will take to compute? What is the largest depth of recursive calls that will be made in Eqn.\ \eqref{eqn:KOH}? This is answered with a brute-force method in {\bf KOHdepth}. {\bf KOHcalls} returns the total number of recursive calls made. But is there an explicit answer?
	
	We will use $n',k'$ to denote the new $n$, $k$ in the recursive call for a given partition. We sometimes also treat them as functions of the index: $n'(i)=(k-i)n-2i+2\sum_{j=0}^{i-1}(i-j)d_{k-j}$. Once a partition only has distinct parts, the recurrence will cease after a final recursive call, since each $k'$ is at most 1. Therefore to maximize depth, we need to maximize using partitions with repeated parts and maximize the size of $n'$, $k'$.
	
	One method is to begin by using $\left[\floor{\frac{k}{2}}^2\right]$ with a $0,1$ as needed and the recursive call with $i=\ceil{\frac{k}{2}}$. Then (since $k'=d_{k-i}=2$) repeatedly use the partition $[1^2]$ and $i=1$. This method calls
	\begin{align*}
		G\left(\left(k-\ceil{\frac{k}{2}}\right)n-2\ceil{\frac{k}{2}}+2\sum_{j=0}^{\ceil{k/2}-1}\left(\ceil{\frac{k}{2}}-j\right)d_{k-j},d_{k-\ceil{k/2}}\right)&=	G\left(\floor{\frac{k}{2}}n-2\ceil{\frac{k}{2}}+2\sum_{j=0}^{\ceil{k/2}-1}0,d_{\floor{k/2}}\right)\\
			&=	G\left(\floor{\frac{k}{2}}n-2\ceil{\frac{k}{2}},2\right).
	\end{align*}
	Then using the partition $[1^2]$, and $i=1$, calls $G(n'-2,2)$. Thus, the total depth of calls is
	\beq
		\ceil{\frac{\floor{k/2}n-2\ceil{k/2}}{2}}+1	=	\ceil{\frac{\floor{k/2}n}{2}}-\ceil{\frac{k}{2}}+1.
	\eeq
	In fact
	\bel
		If $k=1,2$ or if $k\ge4$ is even (odd) and $n\ge2$ ($n\ge4$), then the maximum depth of recursive calls is
		\beq
			\ceil{\frac{\floor{k/2}n}{2}}-\ceil{\frac{k}{2}}+1.
		\eeq
		For $k=3$, the maximum depth of recursive calls is $\ceil{\frac{n}{4}}$.
		\bep
			Base case: $k=1$. There are no recursive calls:
			\beq
				\ceil{\frac{\floor{1/2}n}{2}}-\ceil{\frac{1}{2}}+1	=	\ceil{\frac{0\cdot n}{2}}-1+1	=	0.
			\eeq
			For $k=2$, the only option that yields recursive calls is $[1^2]$ and $i=1$: $G(n-2,2)$. Thus, we will have
			\beq
				\ceil{\frac{n}{2}}	=	\ceil{\frac{\floor{2/2}n}{2}}-\ceil{\frac{2}{2}}+1
			\eeq
			recursive calls, which matches with expected. For $k=3$, the only recursive call that does not immediately terminate is $\lambda=[1^3]$ and $i=2$. The call is
			\beq
				G\left((3-2)n-2\cdot2+2\sum_{j=0}^{2-1}(2-j)d_{3-j},d_{3-2}\right)	=	G(n-4,3).
			\eeq
			Thus, we will have $\ceil{\frac{n}{4}}$ recursive calls. The difference from other $k$ arises because $\floor{\frac{3}{2}}=1$.\\
			
			First consider even $k$ and $n=2$. By Proposition \ref{pro:MaxDepth1}, the depth is
			\beq
				1	=	\frac{kn}{4}-\frac{k}{2}+1	=	\ceil{\frac{\floor{k/2}n}{2}}-\ceil{\frac{k}{2}}+1.
			\eeq
			
			Now let $k\ge4$ and $n\ge3$ ($n\ge4$ if $k$ is odd). Assume the lemma is true for $(n',k')$ s.t. $k'<k$ or $k'=k$ and $n'<n$.\\
			What if we have a partition that repeats part $m$, $\ell\ge2$ times? The recursive call relevant to $m$ uses $i=k-m$ and is
			\begin{equation*}
				G\left(\left(k-k+m\right)n-2\left(k-m\right)+2\sum_{j=0}^{k-m-1}\left(k-m-j\right)d_{k-j},d_{m}\right)	=	G\left(mn-2k+2m+2\sum_{j=m+1}^{k}\left(j-m\right)d_{j},\ell\right).
			\end{equation*}
			By the induction hypothesis, the depth of this call is (note $ml\le k$)
			\begin{align*}
				1-\ceil{\frac{\ell}{2}}+\ceil{\frac{\floor{\ell/2}\left(mn-2k+2m+2\sum_{j=m+1}^{k}\left(j-m\right)d_{j}\right)}{2}}&\le	1-\ceil{\frac{l}{2}}+\ceil{\frac{\floor{l/2}\left(mn-2k+2m+2\sum_{j=m+1}^{k}jd_{j}\right)}{2}}\\
						&\le	1-\ceil{\frac{\ell}{2}}+\ceil{\frac{\floor{\ell/2}\left(mn-2k+2m+2(k-\ell m)\right)}{2}}\\
						&=	1-\ceil{\frac{\ell}{2}}+\ceil{\frac{\floor{\ell/2}m\left(n+2-2\ell\right)}{2}}\\
						&\le	1-\ceil{\frac{\ell}{2}}+\ceil{\frac{m\ell\left(n+2-2\ell\right)}{4}}\\
						&\le	1-\ceil{\frac{\ell}{2}}+\ceil{\frac{k\left(n+2-2\ell\right)}{4}}.
			\end{align*}
			The upper bound is maximized when $\ell$ is minimized: $\ell=2$. This yields an upper bound of $\ceil{\frac{kn}{4}-\frac{k}{2}}$.\\
			
			We now need to confirm that $\lambda=\left[\floor{\frac{k}{2}}^2\right]$ (with $1$ if $k$ is odd), $i=\ceil{\frac{k}{2}}$ achieves this upper bound. The depth of this call is
			\begin{equation}
				1-\ceil{\frac{2}{2}}+\ceil{\frac{\floor{2/2}\left(\floor{\frac{k}{2}}n-2k+2\floor{\frac{k}{2}}\right)}{2}}	=	\ceil{\frac{\left(\floor{\frac{k}{2}}n-2\ceil{\frac{k}{2}}\right)}{2}}	=	\ceil{\frac{\floor{k/2}n}{2}-\ceil{\frac{k}{2}}}.
				\label{eqn:AchievableDepth}
			\end{equation}
			When $k$ is even, this depth directly matches the upper bound given above. We are then left to verify the depth for odd $k$. We adjust the upper bound above by recognizing that if $\ell$ is even, then $m\ell\le k-1$ and splitting into 2 cases:
			\ben
				\item	$\sum_{j=m+1}^{k}jd_{j}=0$. Then the depth of the call is
					\begin{align*}
						1-\ceil{\frac{\ell}{2}}+\ceil{\frac{\floor{\ell/2}\left(mn-2k+2m+2\sum_{j=m+1}^{k}\left(j-m\right)d_{j}\right)}{2}}&=	1-\frac{\ell}{2}+\ceil{\frac{m\ell(n+2)-2k\ell}{4}}\\
							&\le	1-\frac{\ell}{2}+\ceil{\frac{(k-1)(n+2)-2k\ell}{4}}:
					\end{align*}
					maximized for $\ell=2$ and gives an upper bound $\ceil{\frac{(k-1)(n+2)-4k}{4}}=\ceil{\frac{(k-1)n}{4}-\frac{k+1}{2}}$.
				\item	$\sum_{j=m+1}^{k}jd_{j}>0$. Then $m\ell\le k-(m+1)\le k-2$ and the depth of the call is
					\begin{align*}
\hspace{-7mm}						1-\ceil{\frac{\ell}{2}}+\ceil{\frac{\floor{\ell/2}\left(mn-2k+2m+2\sum_{j=m+1}^{k}\left(j-m\right)d_{j}\right)}{2}}&\le	1-\frac{\ell}{2}+\ceil{\frac{\ell\left(mn-2k+2m+2(k-\ell m)-2m\right)}{4}}\\
							&=	1-\frac{\ell}{2}+\ceil{\frac{m\ell\left(n-2\ell\right)}{4}}\\
							&\le	1-\frac{\ell}{2}+\ceil{\frac{(k-2)\left(n-2\ell\right)}{4}},
					\end{align*}
					which again is maximized for $\ell=2$ and gives an upper bound $\ceil{\frac{(k-2)(n-4)}{4}}=\ceil{\frac{(k-1)n}{4}-\frac{k+1}{2}-\frac{2k+n-10}{4}}\le\ceil{\frac{(k-1)n}{4}-\frac{k+1}{2}}$ since $k\ge4$ and $n\ge2$.\footnote{Since we have that $k$ is odd, $k\ge5$ and therefore $n\ge4$; the bound becomes tighter.}
			\een
			The upper bound from either case matches the achievable depth in Eqn.\ \eqref{eqn:AchievableDepth} for odd $k$.
			
			The final case is if $\ell$ and $k$ are both odd. The earlier estimate in the upper bound of $\floor{\frac{\ell}{2}}\le\frac{\ell}{2}$ can be replaced by $\floor{\frac{\ell}{2}}=\frac{\ell-1}{2}$. Thus, the upper bound is at most
			\beq
				1-\ceil{\frac{\ell}{2}}+\ceil{\frac{m(\ell-1)\left(n+2-2\ell\right)}{4}}	\le	1-\ceil{\frac{\ell}{2}}+\ceil{\frac{(k-1)\left(n+2-2\ell\right)}{4}},
			\eeq
			which is now maximized by $\ell=3$ (since we assume $\ell$ is odd)\footnote{Actually, since we ignore $\ell=3$, we can assume $\ell=5$ giving an even worse upper bound.} giving the upper bound: $\ceil{\frac{(k-1)\left(n-4\right)}{4}}-1$, which is worse\footnote{We only need $k\ge1$.} than the upper bound when $\ell=2$. I.e., no recursive call with an odd $\ell$ and $k$ can have larger depth than that listed in Eqn.\ \eqref{eqn:AchievableDepth}. Therefore our chosen $\lambda$, and $i$, achieve the greatest depth: $1+\ceil{\frac{\floor{k/2}n}{2}}-\ceil{\frac{k}{2}}$.
			
			To confirm that this recursive depth actually happens, we must also confirm that $G(n',k')\ne0$ for any other $i$. I.e., we must show that $n'\ge0$. For $i<\ceil{\frac{k}{2}}$
			\beq
				n'	\hspace{-0.15em}=\hspace{-0.15em}	(k-i)n-2i+2\sum_{j=0}^{i-1}(i-j)d_{k-j}	\hspace{-0.15em}=	\hspace{-0.15em}	(k-i)n-2i	\hspace{-0.1em}>\hspace{-0.3em}	\left(k-\ceil{\frac{k}{2}}\right)n-2\ceil{\frac{k}{2}}	\hspace{-0.15em}=\hspace{-0.15em}	\floor{\frac{k}{2}}n-2\ceil{\frac{k}{2}}.
			\eeq
			Since $n\ge3$ ($n\ge4$) when $k$ is even (odd), we have that $n'\ge0$. For $i>\ceil{\frac{k}{2}}$,
			\begin{align*}
				n'	&=	(k-i)n-2i+2\sum_{j=0}^{i-1}(i-j)d_{k-j}	=	(k-i)n-2i+2\left(i-\ceil{\frac{k}{2}}\right)2\\
					&=	(k-i)n+2i-4\ceil{\frac{k}{2}}	=	kn-i(n-2)-4\ceil{\frac{k}{2}}	\ge	kn-(k-1)(n-2)-4\ceil{\frac{k}{2}}\\
					&=	kn-kn+n+2k-2-4\ceil{\frac{k}{2}}	=	n-2+2k-4\ceil{\frac{k}{2}},
			\end{align*}
			which when $k$ is even, $n'\ge0$ since $n\ge3$. And when $k$ is odd, $n'\ge0$, since $n\ge4$.\\
			By induction, the claim holds for all $n,k$.\\
			
			We have glossed over a few details. We assumed that $k'=\ell\ne3$. If $\ell=3$, then this would lead to a different (shallower) recursive call than $\ell=2$:
			\begin{align*}
				\ceil{\frac{n'}{4}}	\le	\ceil{\frac{n'}{2}}	=	\ceil{\frac{\floor{2/2}n'}{2}}-\ceil{\frac{2}{2}}+1.
			\end{align*}
			We also need to ensure that if $k'\ge4$, then $n'\ge3$ if $k'$ is even or $n'\ge4$ if $k'$ is odd. Actually, if $k'\ge4$ is even (odd) and $n'<3$ ($n'<4$) then by Proposition \ref{pro:MaxDepth1} their depth is 1. And since the extra depth of using $\lambda=\left[\floor{\frac{k}{2}}^2,\ k\ \modu{2}\right]$ is\footnote{Recall that $k\ge4$ and $n\ge3$ ($n\ge4$). This is the only spot that needs $n\ge3$ when $k$ is even. Otherwise $n\ge2$ is sufficient.}
			\beq
				\ceil{\frac{\floor{k/2}n}{2}}-\ceil{\frac{k}{2}}	\ge	1,
			\eeq
			$\ell=2$ produces a recursive call that is at least as deep.
		\eep
		\label{lem:MaxDepth}
	\eel
	
	\begin{proposition}
		If $n=1$ or $2$, the maximum depth of recursive calls to $G(n,k)$ is $1$. If $k\ge3$ is odd, then for $n=3$, the maximum depth is $1$ as well.
		\bep
			Consider a partition that repeats part $m$, $\ell\ge2$ times. First consider $n\le2$. The recursive call with $i=k-m$ is
			\begin{align*}
				n'	&=	mn-2k+2m+2\sum_{j=m+1}^{k}\left(j-m\right)d_{j}	\le	4m-2k+2\sum_{j=m+1}^{k}\left(j-m\right)d_{j}\\
					&\le	4m-2k+2(k-m\ell)	=	m(4-2\ell)	\le	m(4-4)	=	0.
			\end{align*}
			Thus, the chain terminates after at most 1 step. For all other partitions (distinct), we already know the chain terminates after 1 step.\\
			Now consider odd $k$ and $n=3$. Note that $d_k=0$ otherwise we have a distinct partition. The recursive call with $i=k-1$ is
			\begin{align*}
				n'	&=	(k-(k-1))3-2(k-1)+2\sum_{j=0}^{(k-1)-1}\left((k-1)-j\right)d_{k-j}\\
					&=	5-2k+2\sum_{j=0}^{k-1}\left(k-j\right)d_{k-j}-2\sum_{j=0}^{k-1}d_{k-j}	=	5-2k+2k-2\sum_{j=0}^{k-1}d_{k-j}	=	5-2\sum_{j=1}^kd_j	=	5-2|\lambda|.
			\end{align*}
			If $|\lambda|\ge3$, $n'<0$. Thus, the only partitions that actually create recursive calls are those with $|\lambda|\le2$. And since $k$ is odd, size 2 partitions will be distinct and therefore terminate after at most 1 step.
		\eep
		\label{pro:MaxDepth1}
	\end{proposition}
	
	The explicit depth formula is {\sc implemented as {\bf KOHdepthFAST(n,k)}}. It is somewhat odd that the depth is not symmetric in $n,k$. $G(n,k)$ is symmetric so we can choose to calculate $G(k,n)$ if $n>k$ and have a shallower depth of recursive calls. The total number of calls may be the same; this was not analyzed.
	
\end{section}

\begin{section}{Altered Recurrence}

	\begin{subsection}{Restricted Partitions}
		We can alter the recurrence of Eqn.\ \eqref{eqn:KOH} in several ways to create $G'(n,k)$ and maintain the property that $G'(n,k)$ is symmetric and unimodal of darga $nk$.\footnote{Dependent on restrictions, one may achieve $G'(n,k)=0$.} One simple way to do this is to restrict the partitions over which we sum. This simply reduces the use of Proposition \ref{pro:SumSymUni}. We can restrict the minimum/maximum size of a partition, min/max integer in a partition, distinct parts, modulo classes, etc.
	
		Suppose we restrict to partitions with size $\le s$ and denote this new function as
		\begin{align*}
			G_s(n,k)	&=	\sum_{(d_1,\ldots,d_k);\sum_{i=1}^k id_i=k;\sum_{i=1}^k d_i\le s} q^{k\left(\sum_{i=1}^kd_i\right)-k-\sum_{1\le j<i\le k}(i-j)d_id_j}	\prod_{i=0}^{k-1}G_s\left((k-i)n-2i+2\sum_{j=0}^{i-1}(i-j)d_{k-j},d_{k-i}\right)\\
					&=	\sum_{(d_1,\ldots,d_k);\sum_{i=1}^k id_i=k;\sum_{i=1}^k d_i\le s} q^{k\left(\sum_{i=1}^kd_i\right)-k-\sum_{1\le j<i\le k}(i-j)d_id_j}\prod_{i=0}^{k-1}G\left((k-i)n-2i+2\sum_{j=0}^{i-1}(i-j)d_{k-j},d_{k-i}\right).
		\end{align*}
		For $k\le s$, $G_s(n,k)=G(n,k)$ since ALL partitions of $k$ will have size $\le s$. And since $d_{k-i}\le s$ by definition of the sum, we can replace $G_s$ by $G$ (the $q$-binomial of Eqn. \eqref{eqn:Qbinomial}) in the product. I found conjectured recurrence relations of order 1 and degree 1 in both $n,k$ for $s=1,2,3,4$. I then conjectured for general $s$ that
	
		\begin{conjecture}
			For $G_s(n,k)$ as described above, we have the recurrence relation
			\begin{equation}
				q^n(q^{k+s}-1)G_s(n+1,k)-q^{k+1}(q^n-q^{s-2})G_s(n,k+1)+(q^n-q^{k+s-1})G_s(n+1,k+1)	=	0.
				\label{eqn:RecurrenceSizeS}
			\end{equation}
			\label{con:sizeS}
		\end{conjecture}
		Eqn.\ \eqref{eqn:RecurrenceSizeS} was verified for $n,k\le20$ and $s\le10$. The bounds were chosen to make the verification run in a short time ($\approx2$ minutes). In the algebra of formal power series, taking the limit as $s\to\infty$ (thus obtaining the original $G(n,k)$ with no restriction), recovers
		\begin{align*}
			q^n(-1)G_\infty(n+1,k)-q^{k+1}(q^n)G_\infty(n,k+1)+(q^n)G_\infty(n+1,k+1)	&=	0,\\
			G(n+1,k)+q^{k+1}G(n,k+1)	&=	G(n+1,k+1).
		\end{align*}
		which matches the recurrence relation in Eqn.\ \eqref{eqn:Recurrence}. Recall that $G_s(n,k)=G(n,k)$ for $k\le s$. It is somewhat surprising that $G(n,k)$ follows many recurrence relations for bounded $k$.
		
		We can try to enumerate all $G_s(n,k)$ using a translated Eqn.\ \eqref{eqn:RecurrenceSizeS}:
		\beq
			G_s(n,k)	=	\frac{1}{q^{n}-q^{k+s-1}}\left(q^{n}(q^{k+s-1}-1)G_s(n,k-1)-q^{k}(q^{n}-q^{s-1})G_s(n-1,k)\right),
		\eeq
		{\sc implemented as \bf{Gs(s,n,k,q)}}. However, this recurrence by itself cannot enumerate all $G_s(n,k)$ as we will eventually hit the singularity-causing $n=k+s-1$ if $n\ge s$. However, we do obtain an interesting relation along that line:
		\beq
			G_s(k+s,k)	=	\frac{1-q^{k+1}}{1-q^{k+s}}G_s(k+s-1,k+1) \hspace{0.8em}\mbox{or}\hspace{0.8em}	G_{n-k+1}(n+1,k)	=	\frac{1-q^{k+1}}{1-q^{n+1}}G_{n-k+1}(n,k+1).
		\eeq
		So far, Eqn.\ \eqref{eqn:RecurrenceSizeS} has been verified for $s=1,2,3$ using the explicit formulas in Eqns. \eqref{eqn:ExplicitSize1}, \eqref{eqn:ExplicitSize2}, and \eqref{eqn:ExplicitSize3}, respectively.
		
		It is quite remarkable that the generalized recurrence in Eqn.\ \eqref{eqn:RecurrenceSizeS} yields polynomials using the initial conditions of Eqn.\ \eqref{eqn:InitialConditions},\footnote{$G_s(n,0)=G_s(0,k)=1$.} let alone unimodal polynomials. Are there other sets of initial conditions (potentially for chosen values of $s$) that will also yield unimodal polynomials when iterated in Eqn.\ \eqref{eqn:RecurrenceSizeS}?\footnote{At least until the singularity.}\\
		
		There are a couple of avenues for attack of Conjecture \ref{con:sizeS}. A good starting place may be to tackle the (easier?) question of whether the recurrence even produces symmetric polynomials.
		
		Since we have base case verification, we could try to use induction to prove the conjecture. This would likely involve writing $G_{s+1}=G_s+g_{s+1}$, where $g_{s+1}$ is the contribution from all of the partitions of size $s+1$:
		\beq
			g_{s+1}(n,k)	=	\sum_{(d_1,\ldots,d_k);\sum_{i=1}^k id_i=k;\sum_{i=1}^k d_i=s+1} q^{ks-\sum_{1\le j<i\le k}(i-j)d_id_j}\prod_{i=0}^{k-1}G\left((k-i)n-2i+2\sum_{j=0}^{i-1}(i-j)d_{k-j},d_{k-i}\right).
		\eeq
		Another possibility, since we have an explicit expression for $G_s$ (and for $G$ if needed), may be to simply plug in to Eqn.\ \eqref{eqn:RecurrenceSizeS} and simplify cleverly. These avenues were not addressed in this paper.\\
		
		A more general question would be to analyze recurrences of type $$a(n,k)F(n+1,k+1)=b(n,k)F(n+1,k)+c(n,k)F(n,k+1)$$ that produce symmetric and unimodal polynomials. What are the requirements on $a,b,c$? If $a_1,b_1,c_1$ and $a_2,b_2,c_2$ both produce unimodal polynomials, what, if anything, can be said about $a_1+a_2$, $b_1+b_2$, $c_1+c_2$? About $a_1a_2$, $b_1b_2$, $c_1c_2$? If one can build from previously known valid recurrences, one could potentially build up to Eqn.\ \eqref{eqn:RecurrenceSizeS} from Eqn.\ \eqref{eqn:Recurrence}.			
		
		\begin{subsubsection}{\tbl Natural Partitions\tbr}
		
			\bel
				Suppose we restrict partitions to only use integers $\le p$. If $p<\frac{2k}{n+2}$ then $G'(n,k)=0$.\footnote{The bound $n+2<2\ceil{\frac{k}{p}}$ is sufficient to make $G'(n,k)=0$.}
				\bep
					We have $d_i=0$ for $i>p$. In the product of Eqn.\ \eqref{eqn:KOH}, consider $i=k-1$ for any restricted partition $\lambda$. Then for the new $n$ in the recursive call we obtain
					\begin{align*}
						n'(k-1)	&=	(k-(k-1))n-2(k-1)+2\sum_{j=0}^{(k-1)-1}((k-1)-j)d_{k-j}\\
								&=	n-2(k-1)+2\sum_{j=k-p}^{k-2}(k-1-j)d_{k-j}\\
																	&=	n-2(k-1)+2\sum_{j=2}^{p}(k-1-(k-j))d_{k-(k-j)}\\
																	&=	n-2(k-1)+2\sum_{j=1}^{p}(j-1)d_{j}\\
																	&=	n-2(k-1)+2k-2|\lambda|\\
																	&\le	n+2-2\ceil{\frac{k}{p}}\\
																	&\le	n+2-2\frac{k}{p}\\
																	&<	n+2-2\frac{k}{2k/(n+2)}	=	0.
					\end{align*}
					Then $G'(n',d_{1})=0$ for any partition and therefore $G'(n,k)=\sum_{\lambda}0=0$.
				\eep
				\label{lem:BoundedPartitions}
			\eel
			What do other natural partitions yield for a recursive call? Let us first consider $\lambda=[1^k]$. Then $d_1=k$ and $d_{i\ne1}=0$. The recursive call is then
			\begin{align*}
				q^{k\left(\sum_{i=1}^kd_i\right)-k-\sum_{1\le j<i\le k}(i-j)d_id_j}		&\prod_{i=0}^{k-1}G\left((k-i)n-2i+2\sum_{j=0}^{i-1}(i-j)d_{k-j},d_{k-i}\right)\\
				&\tab=	q^{k\cdot k-k}G\left(n-2(k-1),k\right)\prod_{i=0}^{k-2}G\left((k-i)n-2i+2\sum_{j=0}^{i-1}(i-j)\cdot0,0\right)\\
				&\tab=	q^{k\cdot(k-1)}G\left(n-2(k-1),k\right)\prod_{i=0}^{k-2}G\left((k-i)n-2i,0\right)	=	q^{k(k-1)}G\left(n-2(k-1),k\right),
			\end{align*}
			since $(k-i)n-2i\ge2n-2(k-2)$ and if $2n-2(k-2)<0$, then $n-2(k-1)<0$. If we define $G'(n,k)$ as only being over partitions equal to $[1^k]$, then $G'(n,k)=q^{nk/2}$ if $2(k-1)|n$ and $G'(n,k)=0$ otherwise. One can confirm the monomial by either looking at how many times $q^{k(k-1)}$ is multiplied or by recognizing that $G'(n,k)$ must be a monomial of darga $nk$.
			
			Now consider $\lambda=[k]$. Then $d_k=1$, $d_{i\ne k}=0$ and the recursive call is
			\begin{align*}
				q^{k\left(\sum_{i=1}^kd_i\right)-k-\sum_{1\le j<i\le k}(i-j)d_id_j}	&\prod_{i=0}^{k-1}G\left((k-i)n-2i+2\sum_{j=0}^{i-1}(i-j)d_{k-j},d_{k-i}\right)\\
				&\tab=	q^{k-k}G\left(kn,1\right)\prod_{i=1}^{k-1}G\left((k-i)n-2i+2\sum_{j=1}^{i-1}(i-j)d_{k-j}+2(i-0)d_k,0\right)\\
				&\tab=	G\left(kn,1\right)\prod_{i=1}^{k-1}G\left((k-i)n,0\right)	=	G(nk,1)	=	\frac{1-q^{nk+1}}{1-q},
			\end{align*}
			since $(k-i)n\ge n\ge0$. So $[k]$ provides the \tbl base\tbr\ of $G(n,k)$: $1+q+\cdots+q^{nk}$.
			
			If we restrict to partitions of size 1, then
			\begin{equation}
				G_1(n,k)	=	\frac{1-q^{nk+1}}{1-q}.
				\label{eqn:ExplicitSize1}
			\end{equation}
			We can use this to verify the recurrence relation in Eqn.\ \eqref{eqn:RecurrenceSizeS} for $s=1$.\\
			
			\bel
				Let us restrict partitions to those with distinct parts. Fix $k\in\N$. Let $\ell=\floor{\frac{1}{2}(\sqrt{8k+1}-1)}$.\footnote{I.e., the index of the largest triangular number $\le k$ or the maximum size of a distinct partition of $k$.} Then $G'(n,k)$ is a polynomial of degree $k$ in $q^n$ for $n\ge2\ell-2$.\footnote{It appears that $n=2\ell-3$ also fits the polynomial, though this was left unproven.}
				\bep
					Since $\lambda=[k]$ is distinct, $G'(n,k)$ is of degree $k$ in $q^{n}$. Let $\lambda$ be a partition with distinct parts. Then $d_{i}\le1$ and $|\lambda|=\sum_{j=1}^kd_k\le l$. And for each recursive call, $G'(n',d_i)=0,1,$ or $\frac{1-q^{n'+1}}{1-q}$. Since $k$ is fixed, it has a fixed number of partitions into distinct parts; thus, it remains to show $G'(n',d_i)\ne0$ by showing $n'(i)\ge0$ for all $i$. We do this by showing $n'(i)$ is monotonic:
					\begin{align*}
						n'(i)			&=	(k-i)n-2i+2\sum_{j=0}^{i-1}(i-j)d_{k-j},\\
						n'(i)-n'(i-1)	&=	2\sum_{j=0}^{i-1}(i-j)d_{k-j}-2\sum_{j=0}^{i-1}(i-j)d_{k-j}+2\sum_{j=0}^{i-2}d_{k-j}+2d_{k-i+1}-n-2\\
										&=	2\sum_{j=0}^{i-1} d_{k-j}-n-2	\le	2\ell-n-2	\le	0.
					\end{align*}
					Thus, the minimum occurs at the largest value of $i$: $i=k-1$, and has value
					\begin{align*}
						n'(k-1)	&=	(k-(k-1))n-2(k-1)+2\sum_{j=0}^{k-1-1}(k-1-j)d_{k-j}\\
								&=	n-2(k-1)+2\sum_{j=0}^{k-1}(k-1-j)d_{k-j}\\
								&=	n-2(k-1)+2\sum_{j=0}^{k-1}(k-j)d_{k-j}-2\sum_{j=0}^{k-1}d_{k-j}\\
								&=	n-2(k-1)+2k-2\sum_{j=0}^{k-1}d_{k-j}	=	n+2-2\ell	\ge	0.
					\end{align*}
				\eep
			\eel
			
			Let us return to the minimum of $n'$ for general partitions. It occurs either at $i=k-1$, or if $n\le2|\lambda|-2-2(\#$ of 1s in $\lambda)$, then at the final $i$ s.t. $n+2\ge2\sum_{j=0}^{i-1} d_{k-j}$: see the simplification of $n'(i)-n'(i-1)$ above. The second difference is $n'(i)-2n'(i-1)+n'(i-2)=2d_{k-i+1}\ge0$. Thus, the extrema of $n'$ found is in fact a minimum. And if we fix the $i$ such that there exists a minimum (or equivalent to the minimum), we obtain a value of
			\begin{align*}
				(k-i)n-2i+2\sum_{j=0}^{i-1}(i-j)d_{k-j}	&\approx	(k-i)\left(2\sum_{j=0}^{i-1} d_{k-j}-2\right)-2i+2\sum_{j=0}^{i-1}(i-j)d_{k-j}\\
											&=	2\sum_{j=0}^{i-1}(k-j)d_{k-j}-2k	\le	2k-2k	=	0.
			\end{align*}
			To avoid the approximation step, consider the next $i$:
			\begin{align*}
				(k-(i+1))n-&2(i+1)+2\sum_{j=0}^{i}(i+1-j)d_{k-j}<	(k-(i+1))\left(2\sum_{j=0}^{i} d_{k-j}-2\right)-2(i+1)+2\sum_{j=0}^{i}(i+1-j)d_{k-j}\\
													&=	2\sum_{j=0}^{i}(k-j) d_{k-j}-2k	\le	2k-2k	=	0.
			\end{align*}
			The recursive call is then seen to be 0 for any partition that has $|\lambda|-|\{1s\mbox{ in }\lambda\}|\ge \frac{n+2}{2}$.
			
			Now we consider $\lambda=[k-\ell,\ell]$ for some $1\le \ell<\frac{k}{2}$. The recursive call is
			\begin{align*}
				q^{k\left(\sum_{i=1}^kd_i\right)-k-\sum_{1\le j<i\le k}(i-j)d_id_j}	&\prod_{i=0}^{k-1}G\left((k-i)n-2i+2\sum_{j=0}^{i-1}(i-j)d_{k-j},d_{k-i}\right)\\
				&\tab=	q^{2k-k-((k-\ell)-\ell)}\prod_{i=0}^{\ell-1}G\left((k-i)n-2i+2\sum_{j=0}^{i-1}(i-j)d_{k-j},0\right)\\
				&\tab\tab	G\left((k-\ell)n-2\ell+2\sum_{j=0}^{\ell-1}(\ell-j)d_{k-j},d_{k-\ell}\right)\\
				&\tab\tab	\prod_{i=\ell+1}^{k-\ell-1}G\left((k-i)n-2i+2\sum_{j=0}^{i-1}(i-j)d_{k-j},0\right)\\
				&\tab\tab	G\left((k-(k-\ell))n-2(k-\ell)+2\sum_{j=0}^{k-\ell-1}(k-\ell-j)d_{k-j},d_{\ell}\right)\\
				&\tab\tab	\prod_{i=k-\ell+1}^{k-1}G\left((k-i)n-2i+2\sum_{j=0}^{i-1}(i-j)d_{k-j},0\right)\\
				&\tab=	q^{2\ell}\prod_{i=0}^{\ell-1}G\left((k-i)n-2i,0\right)\\
				&\tab\tab	G\left((k-\ell)n-2\ell,1\right)\\
				&\tab\tab	\prod_{i=\ell+1}^{k-\ell-1}G\left((k-i)n-2\ell,0\right)\\
				&\tab\tab	G\left(\ell n-2\ell,1\right)\\
				&\tab\tab	\prod_{i=k-\ell+1}^{k-1}G\left((k-i)n+2i-2k,0\right),
			\end{align*}
			which to be nonzero requires that
			\begin{align*}
				(k-i)n-2i	&\ge	0,\\
				(k-i)n-2\ell	&\ge	0,\\
				(k-i)n+2i-2k	&\ge	0.
			\end{align*}
			This is true because
			\begin{align*}
				(k-i)n-2i	\ge	(k-\ell+1)n-2(\ell-1)	\ge	(2\ell-\ell+1)2-2(\ell-1)	&\ge	4,\\
				(k-i)n-2\ell	\ge	(k-(k-\ell-1))n-2\ell	=	(\ell+1)n-2\ell	&\ge	2,	\tab	\mbox{and}\\
				(k-i)n+2i-2k	\ge	(k-(k-1))n+2(k-1)-2k	=	n-2	&\ge	0.
			\end{align*}
			We assume $n\ge2$ otherwise we get $0$ in the recursive call. Thus, the recursive call is
			\beq
				q^{2\ell}\frac{1-q^{(k-\ell)n-2\ell+1}}{1-q}\cdot\frac{1-q^{\ell n-2\ell+1}}{1-q}.
			\eeq
			
			For a visual of how the different partitions can contribute to the final polynomial, see Figure \ref{fig:PartitionDivision}.
			\begin{figure}[H]
				\caption{Contributions from Different Partitions to Produce $G(5,5)$.}
				\includegraphics[width=\textwidth]{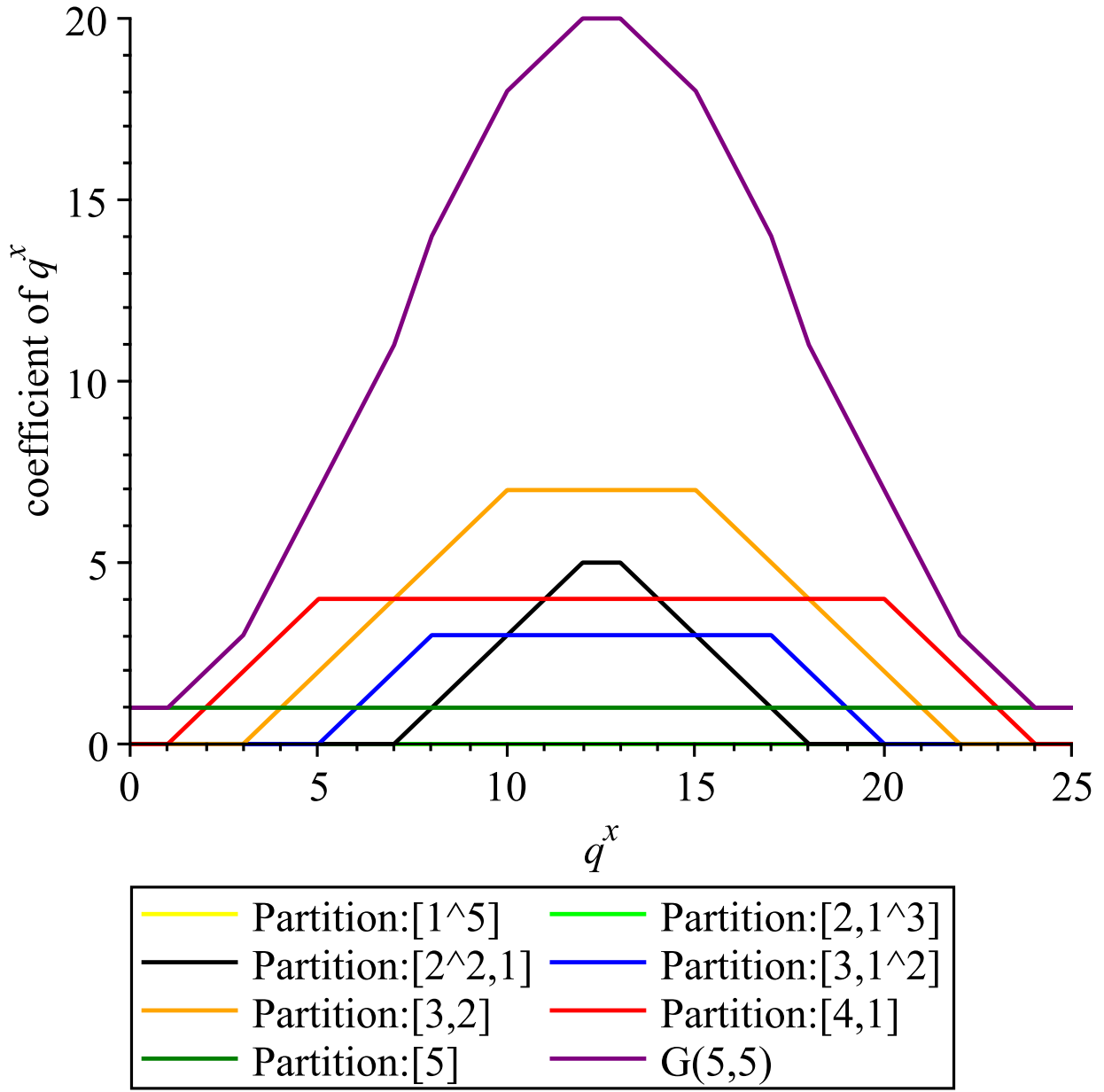}
				 The partitions $[1^5]$ and $[2,1^3]$ both contribute the 0 polynomial.
				\label{fig:PartitionDivision}
			\end{figure}

			\bel
				If we restrict to partitions with consecutive difference $\ge k-d$, and $d<\frac{2}{3}k+1$, then the coefficient of $q^{nk/2}$ in $G'(n,k)$ is $n\frac{D(D+1)}{2}-(D^2-1)$ for $D=\floor{\frac{d}{2}}$.
				\bep
					The bound on $d$ implies that the maximum size of a partition is $2$ by Proposition \ref{pro:RestrictMinDiffSize}. $D$ can be seen to be the relevant number because it is the number of possible size 2 partitions of $k$. $[k],[k-1,1],\ldots,[k-D,D]$ are valid since $k-D-D=k-2\floor{\frac{d}{2}}\ge k-d$. Then
					\begin{equation}
						G'(n,k)	=	\frac{1-q^{nk+1}}{1-q}+\sum_{i=1}^D q^{2i}\frac{1-q^{(k-i)n-2i+1}}{1-q}\cdot\frac{1-q^{in-2i+1}}{1-q}.
						\label{eqn:PartitionsDifferenceKD}
					\end{equation}
					Note that the middle coefficient of $\left(\frac{1-q^{a}}{1-q}\right)\cdot\left(\frac{1-q^{b}}{1-q}\right)$ is going to be $\min(a,b)$. Also note that the darga of the summand in Eqn.\ \eqref{eqn:PartitionsDifferenceKD} is
					\beq
						2\cdot(2i)+[(k-i)n-2i+1-1]+[in-2i+1-1]	=	nk,
					\eeq
					so $q^{nk/2}$ is the middle monomial.\footnote{As one should expect since the summand of Eqn.\ \eqref{eqn:PartitionsDifferenceKD} is a specific case of the summand in Eqn.\ \eqref{eqn:KOH}.} Thus, the coefficient is
					\beq
						1+\sum_{i=1}^D \min\left[(k-i)n-2i+1,in-2i+1\right]	=	1+\sum_{i=1}^D (in-2i+1)	=	\frac{1}{2}\left(  \left( n-2 \right)D +2 \right)  \left( D+1 \right).
					\eeq
				\eep
				\label{lem:CoefficientSize2Partition}
			\eel
			
			\begin{proposition}
				If we restrict to partitions with difference $\ge k-d$ and minimum part $\ell$, then to have a partition of size $s$, it is required that
				\begin{align*}
					d	&\ge	k\frac{(s+1)(s-2)}{s(s-1)}+\frac{2\ell}{s-1},\\
					\ell	&\le	\frac{ds^2-ks^2-ds+ks+2k}{2s},\\
					k	&\le	\frac{(ds-d-2\ell)s}{(s+1)(s-2)}.
				\end{align*}
				\bep
					A partition $\lambda$ of size $s$, minimum part $\ell$, and difference $k-d$ will have total
					\beq
						k	=	\sum_{i=1}^s\lambda_i	\ge	\sum_{i=1}^s \left[\ell+(i-1)(k-d)\right]	=	\frac{s}{2}(2\ell-k+d+(k-d)s).
					\eeq
				\eep
				\label{pro:RestrictMinDiffSize}
			\end{proposition}

		\end{subsubsection}
		
		\begin{subsubsection}{Size 2 Partitions}
		
			Another special partition to consider is $\lambda=\left[\frac{k}{2}^2\right]$ when $k$ is even. Then $d_{k/2}=2$ and $d_i=0$ for $i\ne\frac{k}{2}$ leading to
			\begin{align*}
				q^{k\left(\sum_{i=1}^kd_i\right)-k-\sum_{1\le j<i\le k}(i-j)d_id_j}	&\prod_{i=0}^{k-1}G\left((k-i)n-2i+2\sum_{j=0}^{i-1}(i-j)d_{k-j},d_{k-i}\right)\\
					&=	q^{2k-k-0}\prod_{i=0}^{k/2-1}G\left((k-i)n-2i+2\sum_{j=0}^{i-1}(i-j)d_{k-j},d_{k-i}\right)\\
					&\tab	G\left(\left(k-\frac{k}{2}\right)n-2\frac{k}{2}+2\sum_{j=0}^{k/2-1}\left(\frac{k}{2}-j\right)d_{k-j},d_{k/2}\right)\\
					&\tab	\prod_{i=k/2+1}^{k-1}G\left((k-i)n-2i+2\sum_{j=0}^{i-1}(i-j)d_{k-j},d_{k-i}\right)\\
					&=	q^{k}G\left(\frac{k}{2}n-k,2\right)	\prod_{i=0}^{k/2-1}G\left((k-i)n-2i,0\right)	\prod_{i=k/2+1}^{k-1}G\left((k-i)n-2i+2(i-k/2)2,0\right)\\
					&=	q^{k}G\left(\frac{k}{2}n-k,2\right)	\prod_{i=0}^{k/2-1}G\left((k-i)n-2i,0\right)	\prod_{i=k/2+1}^{k-1}G\left((k-i)n+2i-2k,0\right),
			\end{align*}
			which is nonzero only if $(k-i)n-2i\ge0$ and $(k-i)n+2i-2k\ge0$. By assuming $n\ge2$ (otherwise we get $0$ in the recursive call), we can see that 
			\begin{align*}
				(k-i)n-2i	\ge	(k-(k/2-1))2-2(k/2-1)	=	(k/2+1)2-k+2	&=	4,	\tab\text{and}\\
				(k-i)n+2i-2k	\ge	(k-(k-1))n+2(k-1)-2k	=	n-2	&\ge	0.
			\end{align*}
			Thus, the recursive call is $q^{k}G\left(\frac{k}{2}n-k,2\right)$.\\
			
			While proving Lemma \ref{lem:CoefficientSize2Partition}, Maple produced this useful gem:
			\begin{equation}
				\sum_{i=1}^D q^{2i}\frac{1-q^{(k-i)n-2i+1}}{1-q}\ \cdot\ \frac{1-q^{in-2i+1}}{1-q}	=	q^2\frac{(1+q^{nk-2D})(1-q^{2D})}{(1-q)^2(1-q^2)}-q^{n+1}\frac{(1+q^{nk-nD-n})(1-q^{nD})}{(1-q)^2(1-q^n)}.
				\label{eqn:UsefulGem}
			\end{equation}
			We can then look at using only partitions restricted to maximum size $2$. The possible partitions are $[k]$ and $[k-\ell,\ell]$ for $1\le\ell\le\frac{k}{2}$. We can find the exact form of $G_2(n,k)$ by including all partitions with size $\le2$.
			\bel
				For $n\ge1$,
				\begin{align}
					G_2(n,k)	&=	\bigg(q^{nk+n+4}-q^{nk+n+3}+q^{nk+n+1}-q^{nk+4}-q^{(n-1)k+n+3}+q^{(n-1)k+3}+q^{k+n+1}\nonumber\\
							&\tab	-q^{k+1}-q^n+q^3-q+1\bigg)\bigg/(1-q)^2(1-q^2)(1-q^n).
					\label{eqn:ExplicitSize2}
				\end{align}
				If $n=0$, then $G_2(0,k)=1$.
				\bep
					If $k=$ odd, then, utilizing Eqn.\ \eqref{eqn:UsefulGem},
					\begin{align*}
						G_2(n,k)	&=	\frac{1-q^{nk+1}}{1-q}+\sum_{\ell=1}^{(k-1)/2} q^{2\ell}\frac{1-q^{(k-\ell)n-2\ell+1}}{1-q}\cdot\frac{1-q^{\ell n-2\ell+1}}{1-q}\\
								&=	\frac{1-q^{nk+1}}{1-q}+q^2\frac{(1+q^{nk-2(k-1)/2})(1-q^{2(k-1)/2})}{(1-q)^2(1-q^2)}	-q^{n+1}\frac{(1+q^{nk-n(k-1)/2-n})(1-q^{n(k-1)/2})}{(1-q)^2(1-q^n)}\\
								&=	\frac{1-q^{nk+1}}{1-q}+q^2\frac{(1+q^{nk-k+1})(1-q^{k-1})}{(1-q)^2(1-q^2)}-q^{n+1}\frac{1-q^{n(k-1)}}{(1-q)^2(1-q^n)}.
					\end{align*}
					If $k=$ even, then
					\begin{align*}
						G_2(n,k)	&=	\frac{1-q^{nk+1}}{1-q}+\sum_{\ell=1}^{k/2-1} q^{2\ell}\frac{1-q^{(k-\ell)n-2\ell+1}}{1-q}\cdot\frac{1-q^{\ell n-2\ell+1}}{1-q}+q^{k}G_2\left(\frac{k}{2}n-k,2\right)\\
								&=	\frac{1-q^{nk+1}}{1-q}+q^2\frac{(1+q^{nk-k+2})(1-q^{k-2})}{(1-q)^2(1-q^2)}-q^{n+1}\frac{(1+q^{nk/2})(1-q^{n(k-2)/2})}{(1-q)^2(1-q^n)}\\
								&\tab	+q^k\frac{(1-q^{nk/2-k+1})(1-q^{nk/2-k+2})}{(1-q)(1-q^2)}.
					\end{align*}
					We again utilized Eqn.\ \eqref{eqn:UsefulGem} and noted that $G_2(n',2)=G(n',2)$ and used Eqn.\ \eqref{eqn:Qbinomial}. One can use Maple or any other mathematics software to verify that both cases \tbl simplify\tbr\ to the expression given above.
					
					Notice also that the bound of $n\ge2$ required from characterizing partition calls has been lowered to $n\ge1$. If $n=1$, then the size 2 partition calls will be 0, which is exactly what the extra summands reduce to.
				\eep
				\label{lem:ExplicitSize2}
			\eel
			We can use this to verify the recurrence relation in Eqn.\ \eqref{eqn:RecurrenceSizeS} for $s=2$. To compute $G_2(n,k)$ using the explicit expression, use {\bf G2(n,k,q)}.\\
			
			The expression given in Eqn.\ \eqref{eqn:ExplicitSize2} obfuscates that $G_2(n,k)$ is unimodal of darga $nk$. The expanded expressions in the proof allow for a heuristic argument by reasoning the darga of each summand. This would be a great way to prove unimodality of general functions. However, the useful decomposition is often difficult to ascertain.

		\end{subsubsection}
		
		\begin{subsubsection}{Size 3 Partitions}
		
			By looking at $\lambda=[\ell_1,\ell_2,\ell_3]\vdash k$ for the separate cases $\ell_1=\ell_2=\ell_3$, $\ell_1=\ell_2\ne\ell_3$, $\ell_1\ne\ell_2=\ell_3$, and $\ell_1\ne\ell_2\ne\ell_3$, one can add this to Lemma \ref{lem:ExplicitSize2}, and find an explicit expression for $G_3(n,k)$.\\
			
			We start with $[\ell_1\ne\ell_2\ne\ell_3]$. Then $d_{\ell_i}=1$. We will need $n\ge4$ for nonzero calls as illustrated by
			\begin{align*}
				q^{k\left(\sum_{i=1}^kd_i\right)-k-\sum_{1\le j<i\le k}(i-j)d_id_j}	&\prod_{i=0}^{k-1}G\left((k-i)n-2i+2\sum_{j=0}^{i-1}(i-j)d_{k-j},d_{k-i}\right)\\
					&\tab=	q^{2(\ell_2+2\ell_3)}\prod_{i=0}^{k-\ell_1-1}G\left((k-i)n-2i,0\right)\\
					&\tab\tab	G\left(\ell_1n-2(\ell_2+\ell_3),1\right)\\
					&\tab\tab	\prod_{i=k-\ell_1+1}^{k-\ell_2-1}G\left((k-i)n-2(k-\ell_1),0\right)\\
					&\tab\tab	G\left(\ell_2n-2(\ell_2+\ell_3),1\right)\\
					&\tab\tab	\prod_{i=k-\ell_2+1}^{k-\ell_3-1}G\left((k-i)(n-2)-2\ell_3,0\right)\\
					&\tab\tab	G\left(\ell_3(n-4),1\right)\\
					&\tab\tab	\prod_{i=k-\ell_3+1}^{k-1}G\left((k-i)(n-4),0\right)\\
					&\tab=	q^{2(\ell_2+2\ell_3)}\cdot\frac{1-q^{\ell_1n-2(\ell_2+\ell_3)+1}}{1-q}\cdot\frac{1-q^{\ell_2n-2(\ell_2+\ell_3)+1}}{1-q}\cdot\frac{1-q^{\ell_3(n-4)+1}}{1-q}.
			\end{align*}
			
			Now examine $[k-2\ell,\ell^2]$, for $\ell<\frac{k}{3}$. Then $d_{k-2\ell}=1$, $d_{\ell}=2$. Again, we need $n\ge4$ because the recursive call is
			\begin{align*}
				q^{k\left(\sum_{i=1}^kd_i\right)-k-\sum_{1\le j<i\le k}(i-j)d_id_j}	&\prod_{i=0}^{k-1}G\left((k-i)n-2i+2\sum_{j=0}^{i-1}(i-j)d_{k-j},d_{k-i}\right)\\
					&\tab=	q^{6\ell}\prod_{i=0}^{2\ell-1}G\left((k-i)n-2i,0\right)\\
					&\tab\tab	G\left((k-2\ell)n-4\ell,1\right)\\
					&\tab\tab	\prod_{i=2\ell+1}^{k-\ell-1}G\left((k-i)n-4\ell,0\right)\\
					&\tab\tab	G\left(\ell n-4\ell,2\right)\\
					&\tab\tab	\prod_{i=k-\ell+1}^{k-1}G\left((k-i)(n-4),0\right)\\
					&\tab=	q^{6\ell}\cdot\frac{1-q^{(k-2\ell)n-4\ell+1}}{1-q}\cdot\frac{(1-q^{\ell n-4\ell+1})(1-q^{\ell n-4\ell+2})}{(1-q)(1-q^2)}.
			\end{align*}
			
			We continue by analyzing $[\ell^2,k-2\ell]$, for $\frac{k}{2}>\ell>\frac{k}{3}$. Then $d_{\ell}=2$, $d_{k-2\ell}=1$. $n\ge4$ is needed as usual to show
			\begin{align*}
				q^{k\left(\sum_{i=1}^kd_i\right)-k-\sum_{1\le j<i\le k}(i-j)d_id_j}	&\prod_{i=0}^{k-1}G\left((k-i)n-2i+2\sum_{j=0}^{i-1}(i-j)d_{k-j},d_{k-i}\right)\\
					&\tab=	q^{4k-6\ell}\prod_{i=0}^{k-\ell-1}G\left((k-i)n-2i,0\right)\\
					&\tab\tab	G\left(\ell n-2k+2\ell,2\right)\\
					&\tab\tab	\prod_{i=k-\ell+1}^{2\ell-1}G\left((k-i)n+2i-4(k-\ell),0\right)\\
					&\tab\tab	G\left((k-2\ell)n+8\ell-4k,1\right)\\
					&\tab\tab	\prod_{i=2\ell+1}^{k-1}G\left((k-i)(n-4),0\right)\\
					&\tab=	q^{4k-6\ell}\cdot\frac{(1-q^{\ell n-2k+2\ell+1})(1-q^{\ell n-2k+2\ell+2})}{(1-q)(1-q^2)}\cdot\frac{1-q^{(k-2\ell)n+8\ell-4k+1}}{1-q}.
			\end{align*}
			
			Finally, consider $[\frac{k}{3}^3]$ for $k\cong0\modu{3}$. Then $d_{k/3}=3$. Again, we need $n\ge4$ to discover
			\begin{align*}
				q^{k\left(\sum_{i=1}^kd_i\right)-k-\sum_{1\le j<i\le k}(i-j)d_id_j}	&\prod_{i=0}^{k-1}G\left((k-i)n-2i+2\sum_{j=0}^{i-1}(i-j)d_{k-j},d_{k-i}\right)\\
					&\tab=	q^{2k}\prod_{i=0}^{2k/3-1}G\left((k-i)n-2i,0\right)\\
					&\tab\tab	G\left((n-4)k/3,3\right)\\
					&\tab\tab	\prod_{i=2k/3+1}^{k-1}G\left((k-i)(n-4),0\right)\\
					&\tab=	q^{2k}\frac{(1-q^{(n-4)k/3+1})(1-q^{(n-4)k/3+2})(1-q^{(n-4)k/3+3})}{(1-q)(1-q^2)(1-q^3)}.
			\end{align*}
			
			\bel
				For $n\ge2$, 
				\begin{align}
					G_3(n,k)	&=	\bigg[q^{4k}\cdot{q}^{3} \left( 1-{q}^{n} \right)  \left( q-{q}^{n} \right)
								\tab	-	\tab	q^{3k}\cdot q\left( 1+q \right)  \left( 1-{q}^{n} \right)  \left( q-{q}^{2}+{q}^{5}-{q}^{n} \right)\nonumber\\
								&\hspace{1.2em}-q^{2k}\bigg(q^{nk}\left({q}^{2n}({q}^{9}- {q}^{8}-{q}^{7}+{q}^{6}+ {q}^{5}-{q}^{3}+q)-{q}^{n}({q}^{10}-{q}^{8}+{q}^{6}+{q}^{5})+{q}^{10}\right)\nonumber\\
								&\tab\tab\tab-q^{2n}+{q}^{n}({q}^{5}+{q}^{4}-{q}^{2}+1)-{q}^{9}+{q}^{7}-{q}^{5}-{q}^{4}+{q}^{3}+{q}^{2}-q\bigg)\nonumber\\
								&\hspace{1.2em}+q^k\cdot {q}^{n k}{q}^{3} \left( 1+q \right) \left( 1-{q}^{n} \right)  \left( q^5-q^n+q^{n+3}-{q}^{n+4}\right)
									-	{q}^{n k}{q}^{6} \left( 1-{q}^{n} \right) \left( q-{q}^{n} \right)\bigg]\nonumber\\
								&\tab\bigg/(1-q)^2(1-q^2)^2(1-q^3)(1-q^{n-1})(1-q^{n})q^{2k+1}.
					\label{eqn:ExplicitSize3}
				\end{align}
				For $n<2$, $G_3(n,k)=G_2(n,k)$.
				\bep
					For $n\ge4$,
					\begin{align*}
						G_3(n,k)	&=	G_2(n,k)+\sum_{\ell_3=1}^{\floor{k/3}-1}\sum_{\ell_2=\ell_3+1}^{\floor{(k-\ell_3-1)/2}}q^{2(\ell_2+2\ell_3)}\cdot\frac{1-q^{(k-\ell_3-\ell_2)n-2(\ell_2+\ell_3)+1}}{1-q}\\
								&\hspace{6.95cm}	\cdot\frac{1-q^{\ell_2n-2(\ell_2+\ell_3)+1}}{1-q}\cdot\frac{1-q^{\ell_3(n-4)+1}}{1-q}\\
								&\tab	+	\sum_{\ell=1}^{\floor{(k-1)/3}}q^{6\ell}\cdot\frac{1-q^{(k-2\ell)n-4\ell+1}}{1-q}\cdot\frac{(1-q^{\ell n-4\ell+1})(1-q^{\ell n-4\ell+2})}{(1-q)(1-q^2)}\\
								&\tab	+	\sum_{\ell=\ceil{(k+1)/3}}^{\floor{(k-1)/2}}q^{4k-6\ell}\hspace{-0.2em}\cdot\hspace{-0.2em}\frac{(1-q^{\ell n-2k+2\ell+1})(1-q^{\ell n-2k+2\ell+2})}{(1-q)(1-q^2)}\hspace{-0.2em}\cdot\hspace{-0.2em}\frac{1-q^{(k-2\ell)n+8\ell-4k+1}}{1-q}.
					\end{align*}
					The first nested summation is for the distinct size 3 partitions. The other 2 sums are for partitions of type $[k-2\ell,\ell,\ell]$ and $[\ell,\ell,k-2\ell]$, respectively. If $k\cong0\modu{3}$ then, for the $\left[\frac{k}{3}^3\right]$ partition, we need to add the extra term:
					\beq
						q^{2k}\frac{(1-q^{(n-4)k/3+1})(1-q^{(n-4)k/3+2})(1-q^{(n-4)k/3+3})}{(1-q)(1-q^2)(1-q^3)}.
					\eeq
					Amazingly, Maple is able to simplify the above expression into closed form once one makes assumptions about $k\modu{6}$. It is also necessary to divide the first sum into 2 parts based on the parity of $l_3$. In the divided sum, we are replacing $\ell_3$ by $2\ell_3-1$ and $2\ell_3$, respectively.
					\begin{align*}
						\sum_{\ell_3=1}^{\floor{(\floor{k/3}-1)/2}}\Bigg[&\sum_{\ell_2=2\ell_3}^{\floor{k/2}-\ell_3}q^{2(\ell_2+2(2\ell_3-1))}\cdot\frac{1-q^{(k-(2\ell_3-1)-\ell_2)n-2(\ell_2+2\ell_3)+3}}{1-q}\\
						&\hspace{3.62cm}	\cdot\frac{1-q^{\ell_2 n-2(\ell_2+2\ell_3)+3}}{1-q}\cdot\frac{1-q^{(2\ell_3-1)(n-4)+1}}{1-q}\\
									&\tab+\sum_{\ell_2=2\ell_3+1}^{\floor{(k-1)/2}-\ell_3}q^{2(\ell_2+4\ell_3)}\cdot\frac{1-q^{(k-2\ell_3-\ell_2)n-2(\ell_2+2\ell_3)+1}}{1-q}\\
						&\hspace{4.62cm}	\cdot\frac{1-q^{\ell_2n-2(\ell_2+2\ell_3)+1}}{1-q}\cdot\frac{1-q^{2\ell_3(n-4)+1}}{1-q}\Bigg].
					\end{align*}
					If $k\cong0,1,2\modu{6}$ then $\floor{\frac{k}{3}}-1\cong1\modu{2}$ and so we need to add the final inner sum when $\ell_3=\floor{\frac{k}{3}}-1$:
					\begin{align*}
						\sum_{\ell_2=\floor{k/3}}^{\floor{(k-\floor{k/3})/2}}q^{2(\ell_2+2(\floor{k/3}-1))}&\cdot\frac{1-q^{(k-\ell_2-\floor{k/3}+1)n-2(\ell_2+\floor{k/3})+3}}{1-q}\\
						&\cdot\frac{1-q^{\ell_2n-2(\ell_2+\floor{k/3})+3}}{1-q}\cdot\frac{1-q^{(\floor{k/3}-1)(n-4)+1}}{1-q}.
					\end{align*}
					For each $k\cong0,\ldots,5\modu{6}$, the decomposition simplifies to the same result: Eqn.\ \eqref{eqn:ExplicitSize3}.\footnote{This was only fully done for the $k\cong5\modu{6}$ case. However, the other cases were partly simplified and found to be in agreement empirically; after removing the assumption on $k\modu{6}$, all cases still gave the same result. And the explicit formula in Eqn.\ \eqref{eqn:ExplicitSize3} matches the polynomials found by recursion for all $n=2,\ldots,30$ and $k=0,\ldots30$. If the reader does not wish to take our word for it, they are encouraged to show the other cases for themselves. We highly recommend using Maple or another computer assistant if you do not have time to waste.} 
					
					If $n<4$, then the additions from size 3 partition calls will all be 0 (see the characterization of size 3 partition calls above) and so $G_3(n,k)=G_2(n,k)$. It is a simple exercise to check that the difference evaluates to 0 when $n=2,3$. The difference between the two polynomials is given by
					\begin{align*}
						G_3(n,k)-G_2(n,k)	&=	\bigg[q^{4k}\cdot q^2\left( 1-{q}^{n} \right)  \left( q-{q}^{n} \right)\\
										&\tab	-q^{3k}\cdot q\left( 1+q \right)  \left( 1-{q}^{n} \right)  \left( {q}^{n+3}-{q}^{n+2}-q^n+{q}^{3}\right)\\
										&\tab	+q^{2k}\hspace{-0.13em}\cdot\hspace{-0.13em} q\left( {q}^{3}-{q}^{n} \right)  \left({q}^{nk}\left({q}^{n+1}-{q}^{3}-{q}^{2}+1\right)+q^n(q^3-q-1)+{q}^{2} \right)\\
										&\tab	-q^{k}\cdot{q}^{nk}{q}^{3} \left( 1+q \right)  \left( 1-{q}^{n} \right)  \left(1-q -{q}^{3}+{q}^{n}\right)\\
										&\tab	-{q}^{nk}{q}^{5} \left( 1-{q}^{n} \right)  \left( q-{q}^{n} \right)\bigg]\\
										&\tab\bigg/(1-q)^2(1-q^2)^2(1-q^3)(1-q^{n-1})(1-q^{n})q^{2k}.
					\end{align*}
				\eep
				\label{lem:ExplicitSize3}
			\eel
			We can use this result to verify the recurrence relation in Eqn.\ \eqref{eqn:RecurrenceSizeS} for $s=3$. To compute $G_3(n,k)$ for numeric $n,k$ using the explicit expression, use {\bf G3(n,k,q)}.
		
		\end{subsubsection}
		
		Because we will always have $k'\le s$, it is possible to automate the process of finding $G_s(n,k)$ for any $s$ and have it terminate in finite time. Characterize all of the new size $s$ partitions and add their contributions to $G_{s-1}(n,k)$. However, the number of new types of partitions to consider grows exponentially as $2^{s-1}$; each new part is either the same as the previous part, or smaller. Also, the partitions become more and more complex to describe. One general partition that we can tackle is $\left[\frac{k}{\ell}^\ell\right]$. Then $d_{k/\ell}=\ell$ and the recursive call is
		\begin{align*}
			q^{k\left(\sum_{i=1}^kd_i\right)-k-\sum_{1\le j<i\le k}(i-j)d_id_j}	&	\prod_{i=0}^{k-1}G\left((k-i)n-2i+2\sum_{j=0}^{i-1}(i-j)d_{k-j},d_{k-i}\right)\\
						&=	q^{k(\ell-1)}\prod_{i=0}^{k(\ell-1)/\ell-1}G\left((k-i)n-2i,0\right)\\
						&\tab	G\left(\frac{k}{\ell}n-2k\frac{\ell-1}{\ell},\ell\right)\\
						&\tab	\prod_{i=k(\ell-1)/\ell+1}^{k-1}G\left((k-i)n-2i+2(i-k(\ell-1)/\ell)\ell,0\right)\\
						&	=	q^{k(\ell-1)}G\left(\frac{k}{\ell}(n-2(\ell-1)),\ell\right),
		\end{align*}
		which to be nonzero requires that $n\ge2(\ell-1)$.\\
		
		The final functions found in Lemmas \ref{lem:ExplicitSize2} and \ref{lem:ExplicitSize3} are by no means obviously unimodal. If one instead writes them in decomposed form, they can be reasoned to be unimodal by looking at the darga of each term. Knowing the proper decomposition beforehand allows for easy proof of unimodality, but finding any decomposition is typically extremely difficult. This is the motivation behind reverse engineering.
		
		
		\label{subsec:RestrictedPartitions}
	\end{subsection}

	\begin{subsection}{Adjusted Initial Conditions}
		We can also multiply each recursive call by constant factors or multiply initial conditions while still maintaining base case darga $=nk$. The new initial conditions are then
		\begin{equation}
			G(n<0,k)	=	0,	\hspace{1em}	G(n,k<0)	=	0,	\hspace{1em}	G(0,k)	=	\nu,	\hspace{1em}	G(n,0)	=	\nu,	\tab\text{and}\tab	G(n,1)	=	\nu\frac{1-q^{n+1}}{1-q}.
			\label{eqn:AdjustedInitialConditions}
		\end{equation}
		We can also adjust the recursive call on each $\lambda\vdash k$ by multiplying it by $\rho$. $\rho$ need not be constant; it can be a function of $|\lambda|$, the largest part of $\lambda$, or even more generally, each $\lambda$ gets its own specific value. 
		We choose to \tbl normalize\tbr\ the final answer by dividing by $\nu^k$ (and $\rho$ when it is constant) so that we do not have arbitrary differences. Even after normalizing, there will still be differences from $G(n,k)$.\\
		
		Though we have not yet found meaningful results using the techniques of this section, we are hopeful that in the future we can find a pattern for the difference between the largest and second-largest coefficients.
		\label{subsec:AdjustedInitialConditions}
	\end{subsection}

	\begin{subsection}{Adjusted Recursive Call}
		We return to using the initial conditions in Eqn.\ \eqref{eqn:InitialConditions}. A more complex method is to adjust the recursive call as shown below:
		\begin{align*}
			G(n,k;a,b,c)	&=	\sum_{(d_1,\ldots,d_k);\sum_{i=1}^k id_i=k} q^{k\left(\sum_{i=1}^kd_i\right)-k-\sum_{1\le j<i\le k}(i-j)d_id_j}\\
					&\hspace{3.5cm}	\cdot q^{k a+\left(\sum_{i=1}^kd_i\right)b+c[k(k+1)(n/2-a)-k(k+2b)+\sum_{i=1}^kd_i i^2]}\nonumber\\
					&	\prod_{i=0}^{k-1}G\left((k-i)n-2i+2\sum_{j=0}^{i-1}(i-j)d_{k-j}-2a(k-i)-2b,d_{k-i}-2c;a,b,c\right).
		\end{align*}
		The reduction in darga in the recursive call is exactly balanced by the extra factors of $q$.\\
		However, if $c>0$, since there is some $i$ s.t. $d_{k-i}=0$ for any partition, the product will be $0$ by initial conditions, and thus $G(n,k;a,b,c>0)=0$. And if $c<0$, we obtain infinite recursion for nontrivial choices of $n,k$. Therefore, we should ignore this parameter. The useful adjustment is
		
		\begin{align}
			G(n,k;a,b)	=&	\sum_{(d_1,\ldots,d_k);\sum_{i=1}^k id_i=k} q^{(k+b)\left(\sum_{i=1}^kd_i\right)+k(a-1)-\sum_{1\le j<i\le k}(i-j)d_id_j}\nonumber\\
					&\prod_{i=0}^{k-1}G\left((k-i)(n-2a)-2(i+b)+2\sum_{j=0}^{i-1}(i-j)d_{k-j},d_{k-i};a,b\right),
			\label{eqn:KOHgeneral}
		\end{align}
		{\sc implemented as {\bf KOHgeneral}}.\\
		
		Unless otherwise stated, we are now assuming the type of partition is unrestricted.
		\begin{proposition}
			\ben
				\item	$G(n,k;a,b)$ will have smallest degree $ka+b$.
				\item	If $a+b>\frac{n}{2}$, then $G(n,k;a,b)=0$.
				\item	If $a+b=\frac{n}{2}$, then $G(n,k;a,b)=q^{ka+b}\frac{1-q^{nk+1-2(ka+b)}}{1-q}$.
			\een
			\bep
				\ben
					\item	This follows from inspecting factors of $q$. Note that $k\left(\sum_{i=1}^kd_i\right)-k-\sum_{1\le j<i\le k}(i-j)d_id_j\ge0$ since $G(n,k)$ has smallest degree $1$. Recall it $=0$ for $\lambda=[k]$.
					\item	Let $\lambda\vdash k$. Consider $n'(k-1)$ in the recursive call:
						\begin{align*}
							n'(k-1)	&=	(k-(k-1))(n-2a)-2((k-1)+b)+2\sum_{j=0}^{(k-1)-1}((k-1)-j)d_{k-j}\\
										&=	(1)(n-2a)-2(k-1+b)+2\sum_{j=0}^{k-1}(k-1-j)d_{k-j}\\
										&=	n-2a-2k+2-2b+2\sum_{j=1}^{k}(k-1-(k-j))d_{j}\\
										&=	n-2a-2b+2+\left(-2k+2\sum_{j=1}^{k}jd_{j}\right)-2\sum_{j=1}^kd_{j}\\
										&=	n-2(a+b)+2-2\sum_{j=1}^kd_{j}\\
										&\le	n-2(a+b)+2-2(1)	<	0.
						\end{align*}
						Thus, $G(n'(k-1),d_{1};a,b)=0$ and since $\lambda$ was arbitrary, $G(n,k;a,b)=0$.
					\item	The above has equality only if $|\lambda|=1$, i.e., $\lambda=[k]$. Then
						\begin{align*}
							G(n,k;a,b)	&=	q^{(k+b)\left(\sum_{i=1}^kd_i\right)+k(a-1)-\sum_{1\le j<i\le k}(i-j)d_id_j}\\
									&\tab	\prod_{i=0}^{k-1}G\left((k-i)(n-2a)-2(i+b)+2\sum_{j=0}^{i-1}(i-j)d_{k-j},d_{k-i};a,b\right)\\
									&=	q^{ka+b}G\left((k)(n-2a)-2(b),1;a,b\right)\\
									&\tab	\prod_{i=1}^{k-1}G\left((k-i)(n-2a)-2(i+b)+2\sum_{j=0}^{i-1}(i-j)d_{k-j},0;a,b\right)\\
									&=	q^{ka+b}G\left(k(n-2a)-2b,1;a,b\right)\\
									&\tab	\prod_{i=1}^{k-1}G\left((k-i)(n-2a)-2(i+b)+2(i-0),0;a,b\right).
						\end{align*}
						It only remains to confirm $(k-i)(n-2a)-2b\ge0$, which follows from $a=\frac{n}{2}-b$.
				\een
			\eep
		\end{proposition}
	
	\end{subsection}
	
\end{section}

\begin{section}{OEIS}
	When we take $q\to1^-$ for the normal $q$-binomials, we obtain the common binomial coefficients. What happens if we do that to the modified $G_s$ polynomials? There is no direct combinatorial interpretation (yet), but we can derive some cases to start the search. Recall the explicit forms of $G_1(n,k)$, $G_2(n,k)$, and $G_3(n,k)$ in Eqns. \eqref{eqn:ExplicitSize1}, \eqref{eqn:ExplicitSize2}, and \eqref{eqn:ExplicitSize3}, respectively. Their limits are
	\begin{align*}
		\lim_{q\to1^-}G_1(n,k)	&=	nk+1,\\
		\lim_{q\to1^-}G_2(n,k)	&=	\frac{1}{12}(k+1)(k^2n^2-3k^2n-kn^2+2k^2+9kn-8k+12),\\
		\lim_{q\to1^-}G_3(n,k)	&=	\frac{1}{720}(k+2)(k+1)\bigg(k^3n^3-9k^3n^2-3k^2n^3+26k^3n+42k^2n^2+2kn^3\\
							&\hspace{3.7cm}	-24k^3-153k^2n-33kn^2+162k^2+247kn-378k+360\bigg).
	\end{align*}
	We can also conjecture the form of $G_4(n,n)$ and $G_5(n,n)$ from experimental data:
	\begin{align*}
		\lim_{q\to1^-}G_4(n,n)	&=	\frac{1}{120960}n(n+2)(n+1)\bigg(n^8-32n^7+462n^6-3836n^5+20013n^4\\
							&\hspace{4.75cm}	-66836n^3+140804n^2-171216n+100800\bigg),\\
		\lim_{q\to1^-}G_5(n,n)	&=	\frac{n(n+3)(n+2)(n+1)}{43545600}\bigg(n^{10}-50n^9+1140n^8-15420n^7+136533n^6	-824370n^5\\
							&\hspace{4.5cm}+3436190n^4-9762880n^3+18198936n^2-20242080n+10886400\bigg).
	\end{align*}
	The sequences for $G_1(n,n),\ldots,G_5(n,n)$ are now in the OEIS \cite{A002522,A302612,A302644,A302645,A302646}. Only $s=1$ was already in the OEIS.
	\begin{table}[H]
		\caption{Sequences from the OEIS.}
		\begin{align*}
			\seqnum{A002522}:\hspace{1em}	s=1	\hspace{1em}	&	\tab	1, 2, 5, 10, 17, 26, 37, 50, 65, 82, 101, 122, 145, 170, 197, 226, 257, 290, 325, 362, 401,\\
			\seqnum{A302612}:\hspace{1em}	s=2	\hspace{1em}	&	\tab	1, 2, 6, 20, 65, 186, 462, 1016, 2025, 3730, 6446, 10572, 16601, 25130, 36870, 52656,\\
			\seqnum{A302644}:\hspace{1em}	s=3	\hspace{1em}	&	\tab	1, 2, 6, 20, 70, 252, 896, 2976, 8955, 24310, 60038, 136500, 289016, 575680, 1087920,\\
			\seqnum{A302645}:\hspace{1em}	s=4	\hspace{1em}	&	\tab	1, 2, 6, 20, 70, 252, 924, 3432, 12705, 45430, 152438, 472836, 1352078, 3578680, 8827080,\\
			\seqnum{A302646}:\hspace{1em}	s=5	\hspace{1em}	&	\tab	1, 2, 6, 20, 70, 252, 924, 3432, 12870, 48620, 183755, 683046, 2443168, 8263360, 26184420.
		\end{align*}
	\end{table}
	The sequences will approach the central binomial coefficients \seqnum{A000984} \cite{A000984} since $G_s\to G_\infty=G$.
\end{section}

\begin{section}{Conclusion and Future Work}
	In this work we found a new way to generate unimodal polynomials for which unimodality is, {\it a priori}, very difficult to prove. Many of the recurrences found produce surprisingly unimodal polynomials given the correct initial conditions. The proof is \tbl trivial\tbr\ because of our method of reverse engineering. There are many promising avenues for future work.
	\ben
		\item	
			We were not able to make interesting use of the techniques in Section 5.2. With much more care, it may be possible to glean some interesting results by adjusting initial conditions.
		\item	An ultimate goal would be to find a restriction on the partitions and other parameters that yields non-trivial polynomials that can be expressed as a function in 2 variables. We have a few for the restricted partition size, but more is always better!
		\item	Find other initial conditions that still yield unimodal polynomials when iterated in the recurrence of Eqn.\ \eqref{eqn:RecurrenceSizeS}.
		\item	Try to fully automate the process of finding $G_s(n,k)$ in Section \ref{subsec:RestrictedPartitions}. Use them to verify Eqn.\ \eqref{eqn:RecurrenceSizeS} for larger values of $s$. The ideal result would be to prove Eqn.\ \eqref{eqn:RecurrenceSizeS} for ALL $s$.
		\item	I also tried looking at how restricting to odd partitions compares to restricting to distinct partitions. As these are conjugate sets, I thought there may be a relation in the resulting polynomials. Alas, I could see no connection. One could reexamine these pairs or look at other partition sets that are known to be equinumerous, e.g., odd and distinct partitions compared to self-congruent partitions.
		\item	An appealing result would be to find combinatorial interpretations for any of the restricted polynomials. Particularly the $G_s$ polynomials and the OEIS sequences \seqnum{A302612} \cite{A302612}, \seqnum{A302644} \cite{A302644}, \seqnum{A302645} \cite{A302645}, and \seqnum{A302646} \cite{A302646}.
		\item	Instead of only using an adjusted Eqn.\ \eqref{eqn:KOHgeneral}, one might be able to \tbl arbitrarily\tbr\ combine polynomials of known darga to create another polynomial with known darga. Do this in a similar recursive manner. Instead of using partitions, are there other combinatorial objects that can be applied? What other unimodal and symmetric polynomials are out there to use as starting blocks? What recurrences do they satisfy that can be tweaked to reverse-engineer unimodal polynomials?
		\item Unimodal probability distributions have special properties. One could use these polynomials for that purpose after normalizing by $G'(n,k;q=1)$. Various properties can be discovered using Gauss's inequality \cite{GaussInequality} or others.
		\item	
			Unimodal polynomials can be constructed from multiplication very easily. Is the inverse process of factoring them easier than for general polynomials? If so, there could be applications in cryptography and coding theory where factoring is a common theme.
	\een
	Thank you for reading this paper. I hope you have enjoyed it and can make use of this package.
\end{section}

\begin{section}{Acknowledgements}
	I would like to thank Doron Zeilberger for his direction in this paper. I would also like to thank Cole Franks for his edits and suggestions. And thank you to Michael Saks for his discussion and comments on Conjecture \ref{con:sizeS}.
	
	Thank you to the anonymous referee for their constructive comments to improve the readability and clarity of this paper.
	
	This research was funded by a SMART Scholarship: USD/R\&E (The Under Secretary of Defense-Research and Engineering), National Defense Education Program (NDEP) / BA-1, Basic Research.
\end{section}

\bibliographystyle{amsalpha}
\bibliography{Gnk.bib}



\end{document}